\documentclass[review, 1p]{elsarticle}

\usepackage{mypreamble}
\usepackage{setspace}
\usepackage{tikz}
\usepackage{amssymb}
\usepackage{algorithm}
\usepackage{algpseudocode}
\usepackage{verbatim}  

\usetikzlibrary{shapes, arrows.meta, positioning, calc, shadows}

\definecolor{myblue}{RGB}{0, 114, 178}     
\definecolor{myorange}{RGB}{213, 94, 0}    
\definecolor{mygrey}{RGB}{153, 153, 153}   

\usepackage[nonumberlist,nogroupskip]{glossaries}

\newglossaryentry{7730}
{
name={ISO~7730:2005},
description={ISO~7730:2005 is a thermal comfort standard developed by ISO}
}

\newglossaryentry{55}
{
name={ASHRAE~55--2020},
description={ASHRAE~55--2020 is a thermal comfort standard developed by ANSI and ASHRAE}
}

\makenoidxglossaries

\DeclareMathOperator{\tr}{tr}
\newcommand{\divv}{(\nabla\cdot \mathbf{v})}

\newcommand{\Sv}{\mathbf{S\textsubscript{visc}}}
\newcommand{\Ss}{\mathbf{S\textsubscript{solid}}}
\newcommand{\eps}{\boldsymbol{\varepsilon}}
\newcommand{\bsigma}{\boldsymbol{\sigma}}

\newenvironment{lcases}{\left\lbrace\begin{aligned}}{\end{aligned}\right.}
\journal{Journal of Computational Physics}
\begin{document}

    \begin{frontmatter}

    \title{Meshfree Snow Modelling using a Modified Cam-Clay Approach}

    \author[label1,label2]{Erik Schlesinger\corref{mycorrespondingauthor}}
    \ead{schlesingererik@googlemail.com, erik.schlesinger@volkswagen.de}
    \author[label3]{Chaitanya Sanghavi}
    \author[label3]{Jörg Kuhnert}
    \author[label2]{Carsten Schilde}
    \author[label3,label4]{Pratik Suchde}

    \affiliation[label1]{organization={Volkswagen AG},
             addressline={Berliner Ring 2},
             city={Wolfsburg},
             postcode={38440},
             state={},
             country={Germany}}    
    \affiliation[label2]{organization={Forschungsbereich PVZ Partikelsimulation, Institut für Partikeltechnik, Technische Universität Braunschweig},
             addressline={Franz-Liszt Str. 35A},
             city={Braunschweig},
             postcode={38106},
             state={},
             country={Germany}}
    \affiliation[label3]{organization={Fraunhofer-Institut für Techno- und Wirtschaftsmathematik ITWM},
             addressline={Fraunhofer Platz 1},
             city={Kaiserslautern},
             postcode={67663},
             state={},
             country={Germany}}
   \affiliation[label4]{organization={University of Luxembourg},
             addressline={2 Av. de l'Universite},
             city={Esch-sur-Alzette},
             postcode={4368},
             state={},
             country={Luxembourg}}

    \cortext[mycorrespondingauthor]{Corresponding author}

    \begin{abstract}
    Snow is a complex geomaterial whose macroscopic response is governed by density, temperature, and the topology of its evolving microstructure. Its mechanical behavior spans elastic, plastic, viscous, and failure dominated regimes, imposing significant challenges for numerical methods, which intends to simulate large deformations, evolving free surfaces, and complex boundary interactions. This work presents the first integration of a Modified Cam-Clay constitutive formulation for snow into a purely meshfree strong-form collocation framework based on the Generalized Finite Difference Method. The main methodological contribution is a numerical coupling that combines a global implicit mixed formulation for pressure and velocity with a constitutive return-mapping algorithm. The hydrostatic pressure contribution is obtained from a Poisson equation and subsequently corrected through the Modified Cam-Clay return-mapping procedure, while the deviatoric response is treated semi-implicitly using a numerical viscosity formulation. This partitioned treatment of the volumetric and deviatoric stress contributions enables stable simulations with comparatively large time steps while producing smooth spatial pressure fields. As a result, forces on complex boundary geometries can be evaluated accurately. Numerical results of this coupling illustrate the algorithmic stability of the framework, the effective imposition of boundary conditions, and the suitability of local spatial refinement. The feasibility of applying the framework to vehicle--snow interaction through rigid-body coupling is also illustrated. The presented formulation provides a robust basis for future simulations of dynamic snow loading on vehicle structures.
    \end{abstract}

    \begin{graphicalabstract}
        \includegraphics[width=\linewidth]{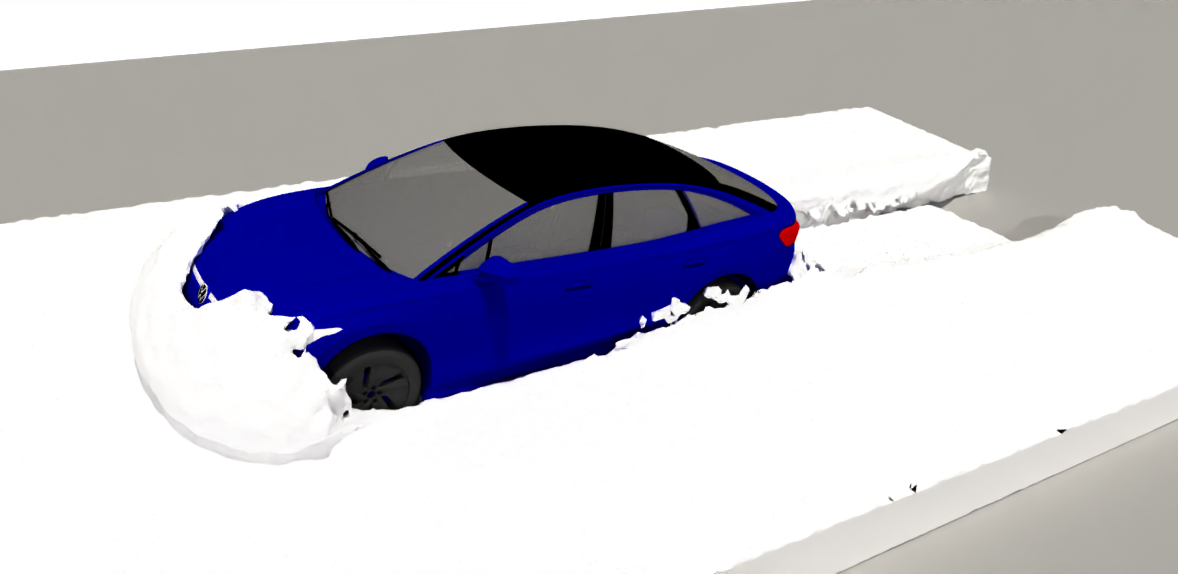}
    \end{graphicalabstract}

    \begin{highlights}
        \item First integration of a Modified Cam-Clay snow model into a purely meshfree strong-form collocation framework. 
        \item A unified implicit pressure-velocity and constitutive return-mapping formulation is developed for snow mechanics. 
        \item Algorithmic stability achieved through semi-implicit deviatoric stress treatment.
        \item Smooth pressure fields are obtained, improving boundary force evaluation.
        \item Complex boundaries, free surfaces, large deformations, and local refinement are handled without remeshing.
        \item The developed framework is coupled with a rigid-body solver for snow-structure interaction simulations.
        \item The formulation provides a basis for vehicle-snow interaction simulations and dynamic load mapping onto structures.
    \end{highlights}


\end{frontmatter}

    \section*{Nomenclature}
\renewcommand{\baselinestretch}{0.75}\normalsize
\renewcommand{\aclabelfont}[1]{\textsc{\acsfont{#1}}}
\begin{acronym}[longest]

    \acro{gfdm}[GFDM]{Generalized Finite Difference Method}
    \acro{fpm}[FPM]{Finite Pointset Method}
    \acro{mcc}[MCC]{Modified Cam Clay Model}
    \acro{csl}[CSL]{Critical State Line}
    \acro{vcl}[VCL]{Virgin Compression Line}
    \acro{urls}[URLs]{Unloading-Reloading Lines}
    \acro{sph}[SPH]{Smooth Particle Hydrodynamics}
    \acro{mpm}[MPM]{Material Point Method}
    \acro{fem}[FEM]{Finite Element Method}
    \acro{pfem}[PFEM]{Particle Finite Element Method}
    \acro{fsi}[FSI]{Fluid Structure Interaction}
    \acro{cog}[COG]{Center Of Gravity}

    \acro{v}[$\mathbf{v}$]{velocity field\acroextra{, m/s}}
    \acro{rho}[$\rho$]{density\acroextra{, kg/m\textsuperscript{3}}}
    \acro{a}[$\mathbf{a}$]{acceleration\acroextra{, m/s\textsuperscript{2}}}

\end{acronym}
\renewcommand{\baselinestretch}{1}\normalsize




    \acresetall

    \section{Introduction}\label{sec:introduction}

Snow is one of the most complex geomaterials occurring in nature. It consists of an ice skeleton and an interstitial pore space, where the diverse morphologies of the ice microstructure significantly influence the overall macroscopic behavior. The vast range of snow types, ranging from highly compressible fresh snow to cohesionless depth hoar and rounded-grain granular snow and even solid ice layers, exhibit a wide array of material properties. These properties are fundamentally governed by bulk density, temperature, and the specific topology of the microstructure. Consequently, developing a representative model remains a significant challenge, requiring a balanced approach between the principles of solid mechanics, soil mechanics and fluid dynamics. For a comprehensive review of the general physical and mechanical properties as well as the modelling of snow, the reader is referred to \cite{Mellor1974, Shapiro1997, Schweizer2003, Blackford2007, Fellin2025}.

In the automotive industry, an accurate representation of boundary conditions, free surfaces, and complex material behavior is essential for evaluating different load cases. These scenarios encompass vehicle interaction with snow, ranging from vehicular snow crossing, structural compression, mechanical loading on various components, as well as snow accumulation. This paper proposes a novel integration of a Modified Cam-Clay constitutive formulation into a strong-form meshfree collocation framework, aimed at capturing snow-specific large deformations, evolving boundaries, and phase transitions. The scope of this research is restricted to a purely mechanical but thermodynamically consistent model that excludes thermodynamic evolution.

From a phenomenological perspective, snow exhibits a dual nature. In the regime of infinitesimal strains and high strain rates, the material behaves as a linear elastic solid, although this elastic domain is notably restricted. Beyond the yield point, the material response transitions into complex behaviors governed by the loading rate. Specifically, rapid loading induces brittle fracture, whereas lower loading rates result in ductile behavior dominated by viscous and plastic dissipation. This pronounced rate dependence and the characteristic ductile-to-brittle transition have motivated the development of a wide array of constitutive models. While early seminal studies approximated snow as a purely elastic medium (e.g., \cite{Smith1972, Heierli2008, Habermann2008}), contemporary frameworks employ advanced constitutive laws to represent its multifaceted response. In particular, elastoplastic formulations are fundamental for describing irreversible compaction and strain hardening or softening \cite{Meschke1996}. To account for the inherent rate-sensitivity of the material, viscoelasticity \cite{Mellor1974, Moos2003, Scapozza2003, Desrues2017} is utilized to capture time-dependent creep, whereas viscoplastic \cite{Cresseri2010} and elasto-viscoplastic formulations \cite{Blatny2024, Vallero2025} describe rate-dependent yielding. Furthermore, the incorporation of anticrack \cite{Gaume2018} and damage mechanics \cite{Moeineddin2024, Moeineddin2026} enables the accurate simulation of structural failure and fracture processes.

Beyond constitutive modeling, selecting a numerically robust simulation method capable of capturing large three-dimensional deformations, evolving free surfaces, solid-to-fluid transitions, and the direct imposition of physical boundary conditions is critical. Traditional grid-based approaches, primarily the \ac{fem}, provide a rigorous mathematical foundation and have been successfully applied to static problems \cite{Meschke1996, Habermann2008, Bartelt2002}. However, \ac{fem} and similar mesh-based techniques are fundamentally limited by the need for continuous remeshing when the material undergoes severe distortion, fragmentation, and topological changes characteristic of snow failure.

These limitations have driven the simulation of snow toward Lagrangian particle methods. Although \ac{sph} has been explored for snow dynamics \cite{Abdelrazek2014, Gissler2020}, it is often hindered by difficulties in accurately enforcing boundary conditions. Similarly, the \ac{pfem} can become computationally demanding due to the necessity of generating a new mesh whenever elements in the Lagrangian background mesh become overly distorted \cite{Li2020, Onate2004, Zorrilla2004}. Consequently, the \ac{mpm} has become one of the leading technique for snow mechanics. It was initially popularized by visual computing applications \cite{Stomakhin2013} and subsequently adapted for physical avalanche mechanics \cite{Gaume2018}. Nevertheless, while \ac{mpm} avoids mesh distortion by using an Eulerian background grid, it can suffer from grid-crossing errors, difficulties in resolving complex, curved boundaries and imposing boundary conditions precisely \cite{Li2020, Liang2024}. Although \ac{mpm} guarantees exact mass conservation through its use of particles with constant mass, this fundamental characteristic inherently introduces computational challenges regarding the insertion and deletion of material points during a simulation \cite{Li2020}.

To address these challenges, this work proposes the integration of a snow constitutive model into a fully meshfree framework based on a strong-form collocation approach, which operates exclusively on a cloud of computational points. To the authors’ knowledge, such a constitutive coupling has not yet been established within a strong-form collocation framework. Specifically, we employ the \ac{gfdm} \cite{benito2001influence, rao2023upwind, zheng2022theoretical}, which has also been referred to as the \ac{fpm} \cite{Drumm2008, Suchde2018}. \ac{gfdm} has already successfully been applied to various industrial problems including water wading or fording of vehicles \cite{Anthony2015}. 

To integrate the constitutive model into a purely meshfree framework, we introduce specific modifications to a numerical scheme originally developed for weakly compressible materials \cite{Kuhnert2014}. By employing a global, implicit mixed formulation for pressure and velocity, the proposed approach yields globally smooth pressure fields evaluated directly at the FPM points. Furthermore, by resolving the hydrostatic component via a Poisson equation, the system ensures an instantaneous equilibrium state. This represents the first application of such a numerical framework to a Cam-Clay model tailored for snow mechanics.

Furthermore, we present what is, to the best of our knowledge, the first \ac{fsi} simulation coupling a Cam-Clay snow model with a rigid-body solver. This test case is explicitly designed to demonstrate the overall algorithmic stability, the efficacy of local spatial refinement, the robust imposition of boundary conditions, and the resulting smooth pressure distributions.

The remainder of this paper is organized as follows. Section~\ref{sec:snowmodel} introduces the governing conservation equations alongside the Cam-Clay constitutive model adapted for snow. Section~\ref{sec:discretization} provides a brief overview of the meshfree \ac{gfdm}. In Section~\ref{sec:MCC_GFDM}, we detail the proposed numerical scheme and the integration of the Cam-Clay model within this framework, concluding with a discussion of the resulting advantages and potential avenues for future research.
Section~\ref{sec:results} is dedicated to numerical results and validation. 

    \section{The Snow Model}\label{sec:snowmodel}

In the context of the present study, we approach snow as a purely mechanical but thermodynamically consistent continuum. 
We adapt the well-established elasto-viscoplastic material model proposed by \cite{Blatny2024} to suit our specific modeling framework. Specifically, we evaluate the elastoplastic Modified Cam Clay model at each collocation point using a return-mapping scheme. The extracted information is used to determine the pressure and velocity in a global monolithic solver. To ensure further numerical stability, we apply a semi-implicit treatment for the deviatoric part of the elastoplastic response, further elaborated in section~\ref{sec:MCC_GFDM}. This formulation yields a smooth, oscillation-free pressure field. The spatial continuity enhances the accuracy and stability of mapping interface forces onto interacting rigid bodies, such as vehicles.

Driven by an evolving yield surface, the elastoplastic formulation models the physical behavior of snow, capturing its hardening response during volumetric compression and strain-softening under both tension and shear-induced dilation. 

We now present the continuous snow model used.

\subsection{The Continuum Equations}\label{subsec:continuum-equations}
The conserved quantities considered are mass, momentum, and total energy. The governing equations formulated in Lagrangian framework reads then 
\begin{align}
    \frac{D\rho}{Dt}&= - \rho \left(\nabla \cdot \mathbf{v}\right) , \nonumber\\
    \frac{D\mathbf{v}}{Dt} &= \frac{1}{\rho}\left(\nabla \cdot \boldsymbol{\sigma}\right) + \mathbf{g},\label{eq:continuum-energy-final} \\
    \frac{De}{Dt} &= \frac{1}{\rho}\left(\boldsymbol{\sigma} : \nabla\mathbf{v}\right). \nonumber
\end{align}
where $D/Dt=\partial/\partial t + (\mathbf{v}\cdot\nabla)$ is the material derivative, $\rho$ is the density, $\mathbf{v}$ the velocity, $e$ is specific energy, $\mathbf{g}$ the acceleration due to external body forces acting on the continuum and $\boldsymbol{\sigma}$ the Cauchy stress tensor. For a comprehensive explanation of how the equations for the primitive variables are derived from the conservative formulation, the reader is referred to \ref{appendix:conservation-eq} and \cite{Nida2005}. For simplification, $\mathbf{g}$ is taken as purely gravitational. Heat sources and heat fluxes are neglected in this work, as we aim to develop a purely mechanical model for snow. However, it should be noted that incorporating these terms in the last equation of system~\eqref{eq:continuum-energy-final} is essential for capturing temperature-dependent phenomena, such as snow melting, which are not considered in the present work. 

To close the system, a \textit{constitutive law} will express the stress tensor $\boldsymbol{\sigma}$ with the specific energy $e$ in terms of the primitive variable $\mathbf{v}$. Towards this end, we introduce a splitting of the Cauchy stress as in \cite{Kuhnert2014} 
\begin{align}
    \boldsymbol{\sigma} &= -p\mathbf{I} + \mathbf{S\textsubscript{visc}} +\mathbf{S\textsubscript{solid}} \label{eq:cauchystress}\\
    &= -(p_{\text{hyd}}+p_{\text{dyn}})\mathbf{I} + \mathbf{S\textsubscript{visc}} +\mathbf{S\textsubscript{solid}},\nonumber
\end{align}
where $\mathbf{I}$ is the identity matrix and we define the pressure $p:=-1/3\tr(\bsigma)$ and the deviator stress $\mathbf{S}:=\bsigma-1/3\tr(\bsigma)$. 

The total pressure $p$ carries all the isotropic bulk terms (mean stresses) of the total stress tensor. In Section~\ref{sec:MCC_GFDM}, we will split the total pressure into a hydrostatic component $p_{\text{hyd}}$ which corresponds to the effects of (effective) body forces and a dynamic component $p_{\text{dyn}}$ which corresponds in theory to dynamic effects but which will also carry numerical corrections to achieve the correct divergence.

The deviator is divided into a viscous deviatoric component $\mathbf{S}\textsubscript{visc}$ and a solid deviatoric component $\mathbf{S}\textsubscript{solid}$ which will obey together with the pressure a hypoelastoplastic law. 

To specify the material model for the Cauchy stress, we must first recall the theory of finite strains. 

\subsection{Finite Strain Theory}\label{subsec:hyperlastic}
As the motivation is to consider snow on a vehicle movement, we also want to consider large deformations. In the theory of finite strains, a deformation is described as a one-to-one mapping from a material coordinate $\mathbf{X}$ in the initial configuration to the spatial coordinate $\mathbf{x}$ in the deformed current configuration. Thus, the components of the deformation gradient $\mathbf{F}$ are the corresponding derivatives $F_{ij}=\partial x_i/\partial X_j$ with $\det(\mathbf{F})>0$. 

From the generalized strain measures in \cite{Seth1964}, we chose the Hencky strain, which is defined with the deformation gradient as $\eps = \frac{1}{2}\log \mathbf{F}\mathbf{F}^T$. The velocity gradient can also be written in terms of the deformation gradient 
\begin{equation}
    \mathbf{L} = \nabla \mathbf{v} = \dot{\mathbf{F}}\mathbf{F}^{-1}.
\end{equation}
Moreover, we define the volumetric strains $\varepsilon_V$ as the trace of the Hencky strain
\begin{equation}
    \varepsilon_V := \tr(\eps) = \log (\det(\mathbf{F})),\label{eq:volumetricstrain}
\end{equation}
and the equivalent deviatoric strains $\varepsilon_D$ as 
\begin{equation}
    \varepsilon_D:= ||\text{dev}\left(\eps\right)||,\label{eq:volumetricstrain}
\end{equation}
where $||\mathbf{A}||:=\sqrt{\mathbf{A}:\mathbf{A}}$ is the Frobenius norm and $\text{dev} (\mathbf{A}):=\mathbf{A}-\tr(\mathbf{A})$ the deviator of an arbitrary second-order symmetric tensor $\mathbf{A}$ \cite{Blatny2024}.

\subsection{Objective Rates}\label{subsec:objective-rate}
The material constitutive equations should be \textit{objective}, which means independent of the frame of reference. To achieve such an objective rate formulation, we use the so-called \textit{logarithmic rate}, which is 
\begin{equation}
    \overset{\circ}{{\mathbf{A}}}^{\log} = \dot{\mathbf{A}} +\mathbf{A}\mathbf{\Omega}^{\log}-\mathbf{\Omega}^{\log}\mathbf{A}.
\end{equation}
for a second-order symmetric tensor $\mathbf{A}$. The logarithmic spin tensor $\mathbf{\Omega}^{\log}$ makes sure that the rate remains objective. In \cite{Xiao1997a, Xiao1997b, Bruhns1999} it is proven that the logarithmic rate of the Hencky strain is the only objective rate that directly coincides with the stretching tensor $\mathbf{D}$, which is defined as the symmetric part of the velocity gradient
\begin{equation}
    \overset{\circ}{{\eps}}^{\log} = \frac{1}{2}\left(\mathbf{L} +\mathbf{L}^T\right) =:\mathbf{D}.
\end{equation}
This also results in the following strain rates for the scalar equivalent strains
\begin{align}
    \dot{\varepsilon_V} &= \tr(\mathbf{D}) = \tr(\mathbf{L})=(\nabla\cdot v),\nonumber\\
    \dot{\varepsilon_D} &=||\text{dev}\left(\dot{\eps}\right)|| = ||\text{dev}\left(\mathbf{D}\right)||.
\end{align}

\subsection{The Hypoelastic Model}\label{subsec:hypoelastic}
For the elastoplastic solid part of the stress tensor, an additive decomposition applies only for infinitesimal strains ($\eps = \eps^E + \eps^P$). This is in general not true for the finite strain theory. There, it is common to decompose the deformation gradient into a gradient of elastic, reversible deformations $\mathbf{F}^E$ and a gradient of plastic, irreversible deformations $\mathbf{F}^P$ in a multiplicative way
\begin{equation}
    \mathbf{F} = \mathbf{F}^E\mathbf{F}^P.\label{eq:deformationdecomposition}
\end{equation}
This work adopts the hyperelastic model proposed by \cite{Blatny2024}, which is formulated using the work-conjugate pair of elastic Hencky strain $\eps^E$ and Kirchhoff stress $\boldsymbol{\tau}:= J\boldsymbol{\sigma}$ (with $J := \det(\mathbf{F})$). This specific formulation offers significant numerical advantages. It naturally circumvents the need for volume corrections during stress integration and avoids complex terms for the scheme of reprojection to the yield surface. It should be noted that alternative approaches utilizing different configurations exist \cite{Simo1998, deSouza2008}, each offering distinct advantages. Accordingly, the hyperelastic constitutive relation is expressed through a strain energy density function
\begin{equation}
    \Psi_{\tau}(\eps^E) := \frac{1}{2}\lambda_{\tau} \tr(\eps^E)^2 + \mu_{\tau} \tr\left((\eps^E)^2\right),\label{eq:Psi-definition}
\end{equation}
where $\lambda_{\tau}$ and $\mu_{\tau}$ are the two Lam\'{e} parameters with respect to the Kirchhoff stress $\boldsymbol{\tau}$. It is noted that the specific internal energy (per unit mass) is related to the strain energy density potential (per unit reference volume) $e=\Psi / \rho_0$ via the initial mass density $\rho_0$. Inserting this strain energy potential into the energy balance (the last equation of system~\eqref{eq:continuum-energy-final}) yields the hyperelastic constitutive law relating the Kirchhoff stress to the Hencky strain, which is provided in detail in \ref{appendix:thermodynamic-consistency}. Since velocity serves as the primary kinematic variable, the objective rate formulation (from section \ref{subsec:objective-rate}) must be employed, which yields a hypoelastic material model
\begin{equation}
    \overset{\circ}{\boldsymbol{\tau}}^{\log} =  \lambda_\tau\tr(\mathbf{D}^E)\mathbf{I} + 2\mu_\tau\mathbf{D}^E.\label{eq:hypoelasticlaw-kirchhoff}
\end{equation}

\subsection{The Modified Cam-Clay Model for Snow}\label{subsec:cam-clay}
We recall equation~\eqref{eq:deformationdecomposition}. Plastic deformations $\mathbf{F}^P$ are assumed to occur when the current stress state reaches the yield surface $f_y=0$, where $f_y$ is the yield function. Consistent with the pull-back to the reference configuration for the material update, the yield surface is formally expressed in terms of two invariants of the Kirchhoff stress, namely the total pressure $p_\tau$ and the second invariant $(J_2)_\tau$, where we recall the Cambridge-notation (consistent with \cite{Blatny2024}) by its square root as   
\begin{equation}
    q_\tau := \sqrt{(J_2)_\tau} = \sqrt{\frac{1}{2}\text{dev}(\boldsymbol{\tau}):\text{dev}(\boldsymbol{\tau})} = \frac{1}{\sqrt{3}}||\boldsymbol{\tau}||_{\text{Mises}}.
\end{equation}
Originally developed for soils, the \textit{Critical State Theory} \cite{Roscoe1968, Schofield1968} postulates that during yielding, a material evolves toward a state where continuous shear deformation occurs without further volumetric changes, the so-called \textit{critical state}. In the context of snow mechanics, the \ac{mcc} model represents the most prominent framework within this category. Its closed yield surface makes it particularly suitable, as it inherently accounts for volumetric failure mechanisms. Building upon this framework, our approach utilizes the yield surface proposed by \cite{Gaume2018}.

\begin{figure}[htb!]
    \centering 
    \begin{tikzpicture}[scale=1.3]
        \draw[thick, ->] (-0.2,0) -- (5.5,0) node[right] {\small $p_\tau$};
        \draw[thick, ->] (1,-1.5) -- (1,3.5) node[above] {\small $q_\tau$};
        
        \draw[ultra thick] (0,0) arc[start angle=180, end angle=0, x radius=2.5cm,
            y radius =2cm];
    
        \node[above] at (3.3,2.8) {\small Hardening};
        \node[above] at (1.75,2.8) {\small Softening};
        \draw[->] (2.3,2.7) -- (1.6,2.7);
        \draw[->] (2.7,2.7) -- (3.4,2.7);
        
        \draw[dashed] (2.5,0) -- (2.5,3.5);
        
        \draw[dotted, ultra thick] (0,0) -- (3,2.4);
        \node at (2.1, 1.3) {\small M};
        \draw[thick] (1.3,1.04) -- (1.9,1.04) -- (1.9,1.52);
        \node at (1.6, 0.85) {\small $1$};
    
        \node[right] at (4.5,1.5) {\small $F(p_\tau,q_\tau, (p_c)_\tau) = 0$};
        
        \node[below] at (0,0) {\small $-\beta (p_c)_\tau$};
        \node[below] at (5,0) {\small $(p_c)_\tau$};
        
        \node[below] at (0.2,-0.8) {\small Tension};
        \node[below] at (2.1,-0.8) {\small Compression};
        \node[below] at (0.2 , 2.4) {\small Shearing};

        \draw[dash dot, ->] (0.8 , 1.8) -- (0.8 , 2.5);
        \draw[dash dot, ->] (0.8,-0.7) -- (0.1,-0.7);
        \draw[dash dot, ->] (1.2,-0.7) -- (1.9,-0.7);
    \end{tikzpicture}
    \caption{The Yieldsurface of the Modified Cam-Clay (\ac{mcc}) Model from \cite{Gaume2018}. Together with a hardening rule, the consolidation pressure $(p_c)_\tau$ regulates as an internal historical variable how the yield surface increases (hardening) or decreases (softening) depending on the load (tension, compression, shearing). $\beta$ is a measure for cohesion, $-\beta (p_c)_\tau$ is the tensile strength and $M$ is the slope of the Critical State Line (thick, dotted), which defines the friction.}
    \label{fig:yieldsurface}
\end{figure}
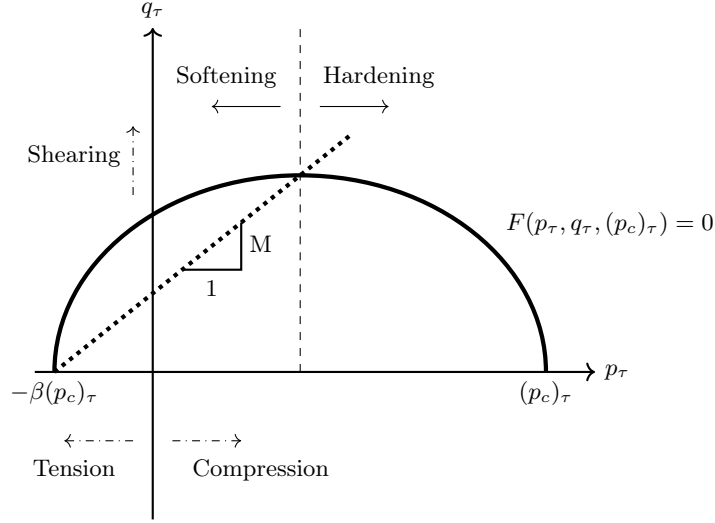

The yield surface consists of all stress states $(p_\tau,q_\tau)$ that satisfy
\begin{equation}
    f_y(\boldsymbol{\tau}) = f_y(p_\tau,q_\tau, (p_c)_\tau) = q_\tau^2 - M^2 (p_\tau+\beta (p_c)_\tau)((p_c)_\tau-p_\tau) \overset{!}{=} 0,
\end{equation}
where the internal historical variable $(p_c)_\tau$ is the so-called consolidation pressure (transformed to the reference configuration). The consolidation pressure represents the critical pressure that initiates volumetric yielding or failure under purely hydrostatic loading. Because this volumetric yielding induces concurrent hardening or softening (dictated by the prevailing stress state) this pressure is defined as an internal variable, governed by a hardening rule based on the plastic volumetric strain. $M$ is the slope of the \ac{csl}. The \ac{csl} corresponds to an evolving Drucker-Prager yield surface $q_\tau(p_\tau)=M(p_\tau+\beta (p_c)_\tau)$, defining the critical state as the apex of the yield function. To account for cohesion, the dimensionless material parameter $\beta$ is introduced. This parameter dictates also the tensile strength of the material, establishing the lower hydrostatic limit of the yield surface at $-\beta (p_c)_\tau$. An illustration of the described yield surface can be found in figure~\ref{fig:yieldsurface}.
The complex architecture of the ice matrix of snow is treated here abstractly through phenomenological internal state variables. These variables govern the initial cohesive strength and the evolution of the yield surface, representing the statistical degradation or densification of the grain bonds, without explicitly resolving the microscopic geometry.
We would like to mention that there are also more sophisticated variants of the \ac{mcc} yieldsurface, e.g. in \cite{Vallero2025}.

The hardening law is based on the typical results of an isotropic compression test \cite{Blatny2024}.  Given that plastic volumetric strains dominate the mechanical response of snow, the elastic contribution can be neglected, allowing for the approximation $\eps_V \approx \eps_V^P$. Furthermore, within an incremental formulation, we may assume that $1 - \exp(\Delta \eps_V) \approx \Delta \eps_V$, since $|\Delta \eps_V| \ll 1$ per time increment. The hardening law results in 
\begin{align}
    (p_c)_\tau = (p_c^{\text{old}})_\tau \exp\left(-\xi \Delta \eps_V^P\right),\label{eq:pc_update}
\end{align}
where the hardening parameter is defined as $\xi:= \rho_{\text{ice}}/(\kappa_P \rho_0)$, where $\rho_{\text{ice}}$ is the density of ice, $\kappa_P$ the compressibility in the plastic regime and $\rho_0$ the initial density. 

The hardening is determined by the amount of $\eps_V^P$. On top of that the amount of plastic volumetric strains is specified with a plastic flow rule. 

The usual convention is to chose the plastic rate of deformation
\begin{equation}
    \mathbf{D}^P = ||\mathbf{D}^P|| \frac{\partial g/\partial\boldsymbol{\tau}}{||\partial g/\partial\boldsymbol{\tau}||}
\end{equation}
in the direction normal to the plastic potential $g$. The typical choice in most of the \ac{mcc} models is a plastic potential $g=f_y$ associated with the yield function $f_y$ \cite{Gaume2018, Blatny2024}. We recall that this is not the case for the non-associated flow rule in \cite{Vallero2025}. Moreover, we define the equivalent plastic strain rate as 
\begin{equation}
    \dot{\gamma} = ||\mathbf{D}^P|| \label{eq:equivplasticstrain}    
\end{equation}
which leads to the associative plastic flow rule 
\begin{align}
    \mathbf{D}^P &= \dot{\gamma} \frac{\partial f_y/\partial\boldsymbol{\tau}}{||\partial f_y/\partial\boldsymbol{\tau}||}.\label{eq:plasticrate}
\end{align}
Equation~\ref{eq:plasticrate} can be splitted into the equivalent plastic shear strain rate
\begin{equation}
    \dot{\eps}_D^P = ||\text{dev}(\mathbf{D}^P)|| = \frac{1}{\sqrt{3}}\dot{\gamma} \frac{\partial f_y/\partial q_\tau}{||\partial f_y/ \partial\boldsymbol{\tau}||}\label{eq:plasticflowrule_dev}
\end{equation}
and equivalent plastic volumetric strain rate
\begin{equation}
     \dot{\eps}_V^P=\tr(\mathbf{D}^P)= - \dot{\gamma} \frac{\partial f_y/\partial p_\tau}{||\partial f_y/ \partial\boldsymbol{\tau}||}. \label{eq:plasticflowrule_vol}
\end{equation}

In summary, if we apply more compression than shear to a sample volume, we are below the \ac{csl}, increasing the consolidation pressure, which makes the yield function grow and the material harder (compare again with figure~\ref{fig:yieldsurface}). At the same time, the volume decreases and the density increases. In contrast, above the critical state line, where we apply tensile forces or strong shearing in relation to compression, the consolidation pressure decreases, the yield function shrinks, and the material softens. At the same time, the volume increases and the density decreases. The critical state is the point at which the current yield surface and \ac{csl} intersect. At this point, neither hardening nor softening occurs, and the material is yielding without any changes in volume. 

\subsection{Viscous Stress Tensor}
To complete the constitutive formulation, we return to the decomposition of the stress tensor presented in equation~\eqref{eq:cauchystress}. It now remains to explicitly specify the deviatoric viscous stress $\Sv$ component to the strain-rate $\overset{\circ}{{\eps}}^{\log}$. This is done with the usual linear, isotropic assumptions (see, for example \cite{Batchelor1967}). We are only interested in the second viscous Lam\'{e} parameter, which is the dynamic viscosity $\eta$ and relates the shear parts. With these assumptions we get 
\begin{equation}
    \Sv = \eta \cdot\text{dev}(\overset{\circ}{{\eps}}^{\log}) =\eta \cdot \text{dev}\left(\mathbf{D}\right).\label{eq:Sv}
\end{equation}
To ensure thermodynamic consistency, the viscous component of the stress tensor governs internal dissipation $\mathcal{D}$, acting as a rate-dependent plastic mechanism. This formulation is conceptually analogous to established overstress models \cite{Blatny2024}, which explicitly permit a temporary excursion of the stress state beyond the static yield surface. Explicit modelling of the viscous component will be part of future work.

    \section{The Generalized Finite Difference Method}\label{sec:discretization}

The numerical framework adopted in this work is a fully Lagrangian Generalized Finite Difference Method (GFDM) \cite{liszka_finite_1980}, which is a meshfree strong-form collocation. In contrast to mesh-based discretizations, the computational domain is represented solely by a cloud of numerical points, allowing the method to naturally accommodate evolving geometries, large deformations, and topological changes without remeshing. This property makes the approach particularly suitable for highly deformable materials such as snow.

\subsection{Domain discretization}\label{subsec:pointcloud}

\begin{figure}[htb!]
\centering
\begin{tikzpicture}[scale=0.8]

\shade[shading=radial, inner color=white, outer color=gray!40] (0,0) circle (3cm);

\coordinate (Center) at (0,0);
\coordinate (InternalJ) at (1.2, -0.8); 

\pgfmathsetseed{42} 

\foreach \x in {-3.5, -2.8, ..., 3.5} {
    \foreach \y in {-3.5, -2.8, ..., 3.5} {
        
        \pgfmathsetmacro{\jx}{\x + rand*0.25}
        \pgfmathsetmacro{\jy}{\y + rand*0.25}
        
        \pgfmathsetmacro{\distSq}{\jx*\jx + \jy*\jy}
        \pgfmathsetmacro{\distJ}{(\jx-1.2)*(\jx-1.2) + (\jy+0.8)*(\jy+0.8)}
        
        \pgfmathsetmacro{\exclArrow}{(\jx > 0.2) && (\jx < 2.9) && (\jy > 0.0) && (\jy < 1.9) && (-0.5*\jx + 0.866*\jy > -0.3) && (-0.5*\jx + 0.866*\jy < 0.7) ? 1 : 0}
        
        \pgfmathsetmacro{\exclXi}{(\jx > -0.6) && (\jx < 0.6) && (\jy > -0.8) && (\jy < -0.1) ? 1 : 0}
        
        \pgfmathsetmacro{\exclXj}{(\jx > 0.8) && (\jx < 1.6) && (\jy > -1.4) && (\jy < -0.8) ? 1 : 0}
        
        \pgfmathparse{\exclArrow + \exclXi + \exclXj > 0 ? 1 : 0}
        \ifnum\pgfmathresult=0
        
            \pgfmathparse{(\jx > -3.8) && (\jx < 3.8) && (\jy > -3.8) && (\jy < 3.8) ? 1 : 0}
            \ifnum\pgfmathresult=1
                
                \pgfmathparse{\distSq > 0.4 ? 1 : 0}
                \ifnum\pgfmathresult=1
                    \pgfmathparse{\distJ > 0.3 ? 1 : 0}
                    \ifnum\pgfmathresult=1
                        
                        \pgfmathparse{\distSq <= 9.1 ? 1 : 0}
                        \ifnum\pgfmathresult=1
                            \fill[myorange] (\jx, \jy) circle (2pt);
                        \else
                            \fill[mygrey] (\jx, \jy) circle (2pt);   
                        \fi
                    \fi
                \fi
            \fi
        \fi
    }
}

\fill[myblue] (Center) circle (2.5pt);
\node[below, text=myblue, yshift=-2pt] (XiLabel) at (Center) {$\mathbf{x}_i$};

\fill[myorange] (InternalJ) circle (2.5pt);
\node [below, text=myorange, xshift=-10pt, yshift=1pt] at (InternalJ) {$\mathbf{x}_j$};

\coordinate (RadiusEnd) at (30:3cm);
\draw [-stealth, thick] (Center) -- (RadiusEnd);
\node [rotate=33, pos=0.5, above, sloped, xshift=35pt, yshift=-2pt] {$h(\mathbf{x}_i, t)$};

\coordinate (LegendStart) at (5, 0.5);

\fill[myblue] (LegendStart) ++(0, 0) circle (2pt) node [right=8pt] {: Collocation point $\mathbf{x}_i$};

\fill[myorange] (LegendStart) ++(0, -1.2) circle (2pt) node [right=8pt, text width=6cm] {: Neighbor points $\mathbf{x}_j$, $j \in S_i$};

\fill[mygrey] (LegendStart) ++(0, -2.4) circle (2pt) node [right=8pt] {: Points outside support domain of $\mathbf{x}_i$};

\end{tikzpicture}
\caption{Illustration of the local support domain and neighbor point selection. The central collocation point $\mathbf{x}_i$ (blue) interacts with its neighbor points $\mathbf{x}_j$ (orange) located within the spherical support domain defined by the smoothing length $h(\mathbf{x}_i, t)$. External points outside this domain do not influence the local approximations.}
\label{fig:gfdm-smoothinglength}
\end{figure}
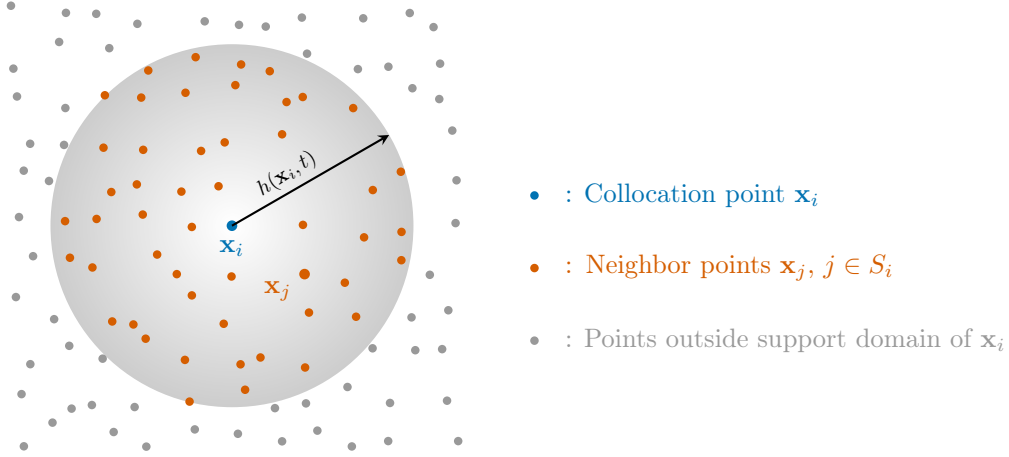

The computational domain is discretized into $N$ computational points, consisting of both interior and boundary points. Their spatial distribution is generated using an advancing front procedure \cite{lohner_advancing_1998, suchde2023point}. The local point density is prescribed through a smoothing length or resolution $h$, which determines the desired spatial refinement, as illustrated in figure~\ref{fig:gfdm-smoothinglength}. In contrast to mass-particle based meshfree methods such as \ac{sph}, the density of the point cloud only indicates the resolution of the discretization and is not related to the physical density \cite{Suchde2023}. 

To maintain a regular point cloud and avoid excessive clustering, inter-point distances are controlled through the parameters $r_{\mathrm{min}}$, $h$, and $r_{\mathrm{max}}$. Specifically, no two points are placed closer than $r_{\mathrm{min}}h$, while at the same time it is ensured that each point has at least one neighbouring point within a sphere of radius $r_{\mathrm{max}}h$. This construction guarantees both numerical stability and sufficient local support for the subsequent approximation of differential operators. Details of the point generation strategy and the enforcement of these distance constraints can be found in \cite{Suchde2019, suchde2023point}.
\subsection{Derivative computation}\label{subsec:gfdm}

Within GFDM, spatial derivatives are approximated directly at each point by constructing local discrete differential operators from neighbouring function values \cite{halada2025overview, Suchde2018}. For a scalar field $f$, the derivative at point $i$ is written as

\begin{equation*}
    \partial^* f (\mathbf{x}_i) \approx \tilde{\partial}^*_i f = \sum_{j \in S_i} c^*_{ij} f_j \,.
\end{equation*}

Here, $S_i$ denotes the local support of point $i$, containing neighbouring points $j$ used for the approximation. The symbol $*$ represents the differential operator under consideration, for example first derivatives in spatial directions or the Laplacian. The notation $\partial^*$ refers to the continuous operator, whereas $\tilde{\partial}^*$ denotes its discrete approximation.

The stencil coefficients $c^*_{ij}$ are obtained from a weighted least-squares procedure, in which consistency conditions are imposed on a local polynomial basis. In this way, the discrete operator reproduces derivatives of selected monomials exactly up to a prescribed order while maintaining numerical robustness on irregular point clouds. We refer to \cite{halada2025overview, Suchde2018, Kuhnert2014} for more details on derivative computation in \ac{gfdm}s.

    \section{Modified Cam-Clay Model in GFDM}
\label{sec:MCC_GFDM}
To establish a robust solver for snow mechanics, the original numerical scheme for weakly compressible media was modified to accommodate the Cam-Clay model, thereby enabling the resolution of large plastic volumetric strains. The adapted framework is detailed in the following section, building upon the foundational methodology originally proposed in \cite{Kuhnert2014}.

\subsection{Numerical Time Discretization Scheme}\label{subsec:scheme}
This section outlines the computational sequence of a single time step. Given the known state variables at the previous step (referred to as $n$), namely density $\rho^n$, velocity $\mathbf{v}^n$, pressure $p^n=p_{\text{hyd}}^n+p_{\text{dyn}}^n$, and the deviatoric elastoplastic stress tensor $\mathbf{S}_{\text{solid}}^n$, a predictor-corrector scheme is employed to compute the updated quantities for the subsequent time increment $\Delta t$.

As it is standard in elastoplastic frameworks, the integration step is initiated by an elastic predictor phase. A preliminary trial pressure is computed under the assumption of a weakly compressible response, yielding an uncorrected trial velocity divergence. For the deviatoric response, $\text{dev}\left(\mathbf{D}\right)$ is utilized to evaluate an elastic trial deviatoric stress. Both the trial pressure and the trial deviatoric stress are subsequently subjected to a plastic return mapping algorithm to enforce yield compliance. The theoretical decomposition of the trial state into independent volumetric and deviatoric components is strictly justified by the assumption of an isotropic hyperelastic response. Governed by the two Lamé parameters, the elastic formulation naturally precludes any cross-coupling between pressure and deviatoric shear during the trial phase.

\subsubsection{Elastic Trial}\label{subsubsec:elastic-trial}
The scheme initiates with the computation of a preliminary hydrostatic pressure field required to balance the body forces, namely gravity, and the divergence of the deviatoric stress. To achieve this, an elastic volumetric compressibility $\kappa_E$ in the Cauchy framework is needed. This bulk moduli can be extracted by transforming the Kirchhoff stress rate in equation~\eqref{eq:hypoelasticlaw-kirchhoff} into the corresponding Cauchy stress rate 
\begin{equation}
    \overset{\circ}{\boldsymbol{\sigma}}^{\log} =  \lambda\tr(\mathbf{D}^E)\mathbf{I} + 2\mu\mathbf{D}^E - \tr(\mathbf{D}^E)\boldsymbol{\sigma}.\label{eq:hypoelasticlaw-cauchy}
\end{equation}
This yields the corresponding elastic compressibility and shear modulus 
\begin{align}
    \kappa_E &= \frac{1}{J K_\tau+p}, \nonumber\\
    \mu_E &= \frac{\mu_\tau}{J}. \label{eq:drhodp}
\end{align}
For a comprehensive explanation of how the moduli are derived from the transformation, the reader is referred to \ref{appendix:thermodynamic-consistency}.
The Poisson-type equation of the hydrostatic pressure then reads
\begin{align}
-\frac{\kappa_E}{\Delta t^2} p_{\text{hyd,E}}^{n+1}+\hat{\nabla}\cdot \left(\frac{{\rm 1}}{\rho } \hat{\nabla} p_{hyd,E}^{n+1} \right)=-\frac{\kappa_E}{\Delta t^2} p_{\text{hyd}}^{n} + \hat{\nabla}\cdot \left(\frac{{\rm 1}}{\rho } \left(\hat{\nabla}\cdot \mathbf{S}_{\text{solid}}^{n} \right)_{}^{T} \right)+ \hat{\nabla}\cdot \mathbf{g}.\label{eq:hydpressure2} 
\end{align}

The preliminary dynamic pressure is assumed as a fraction $C_{\text{dyn}}\in [0,1]$ of the dynamic pressure from the old timestep
\begin{equation}
    \tilde{p}_{\text{dyn,E}}^{n+1} =C_{\text{dyn}}p_{\text{dyn}}^{n}\label{eq:intermediatepressure}
\end{equation}

We linearize the stress update by decomposing the preliminary stress tensor $\tilde{\mathbf{S}}_{\text{solid}}^{n+1}$ into an explicit history-dependent part and an implicit increment
\begin{align}
\hat{\nabla} \cdot \mathbf{S}_{\text{solid}}^{n+1} \approx\hat{\nabla} \cdot \mathbf{S}_{\text{solid}}^n + \hat{\nabla} \cdot \left( 2\mu_{eff} \Delta t \, \mathbf{D}(\mathbf{v}^{n+1}) \right).
\end{align}
The previous stress state $\mathbf{S}_{\text{solid}}^n$ is treated explicitly as a force term on the right hand side, preserving the accumulated plastic history. The increment is approximated by employing an effective shear modulus $\mu_{\text{eff}}$. While this modulus coincides with the standard elastic shear modulus $\mu_E$ within the elastic regime, it is systematically reduced upon reaching the yield surface to account for plastic flow. This introduces an effective (numerical) viscosity 
\begin{align}
    \hat{\eta} :=& \eta +\Delta t \mu_{\text{eff}}\label{eq:effectiveviscosity}
\end{align} 
into the implicit system matrix. Note that we deliberately omit the objective rate rotation terms in the semi-implicit operator. This approximation is justified as it strictly affects only the tangent stiffness matrix of the linear solver, not the constitutive accuracy of the stress state itself, which is updated separately using the full objective rate formulation in the subsequent correction steps.

Governed by the balance of momentum, the intermediate total pressure field $\tilde{p}_E^{n+1}=\tilde{p}_{\text{dyn,E}}^{n+1}+p_{\text{hyd,E}}^{n+1}$ induces a provisional velocity $\tilde{\mathbf{v}}^{n+1}$, which in turn implies a theoretical nominal trial divergence $\overline{\divv}_E$ via mass conservation. We postulate the existence of a correction pressure $c^{n+1}$ designed to enforce a target elastic divergence $\divv_{tar,E}$, thereby strictly satisfying both mass and momentum balances. Consequently, this correction pressure functions as a Lagrange multiplier within a monolithic system of equations formulated for the elastic trial pressure and trial velocity. The system reads
\begin{equation}
    \begin{lcases}
        \frac{\tilde{\mathbf{v}}^{n+1}-\mathbf{v}^{n}}{\Delta t}+ \frac{{\rm 1}}{\rho} \hat{\nabla} c^n+1-\frac{{\rm 1}}{\rho } \left(\hat{\nabla}\cdot \mathbf{S}_{\text{visc}}\left(\mathbf{v}^{n+1}, \hat{\eta}\right) \right)_{}^{T} = &-\frac{{\rm 1}}{\rho } \hat{\nabla} \tilde{p}_E^{n+1}+\frac{{\rm 1}}{\rho } \left(\hat{\nabla}\cdot \mathbf{S}_{\text{solid}}^{n} \right)_{}^{T}+\mathbf{g}\\
        -\frac{\kappa_E}{\Delta t}c^{n+1}+ \hat{\nabla}\cdot\left(\frac{\Delta t_{\text{virt}}}{\rho} \hat{\nabla} c^{n+1}\right)- \hat{\nabla}\cdot\tilde{\mathbf{v}}^{n+1} =&-(\overline{\nabla\cdot\mathbf{v}})_E,&\label{eq:vp-scheme}
    \end{lcases}
\end{equation}
where $\Delta t_{\text{virt}}$ is a virtual timestepsize for subcycling of the correction pressure. For a comprehensive derivation of these equations, the reader is referred to the analogous formulation presented in \cite{Kuhnert2014}.

Having established a smooth, global elastic pressure field, the scheme proceeds to evaluate the elastic deviatoric stress $\mathbf{S}_{\text{solid,E}}^{n+1}$. To ensure material frame indifference, this computation is driven by the velocity field, specifically, the rate of deformation tensor $\text{dev}(\mathbf{D}^{n+1})$, integrated via the objective stress rate of section~\ref{subsec:hypoelastic} 
\begin{equation}
    \mathbf{S}_{\text{solid,E}}^{n+1} =\mathbf{S}_{\text{solid}}^{n} + \Delta t \left(2 \mu_E \text{dev}\left(\mathbf{D}^{n+1}\right) - \mathbf{S}_{\text{solid}}^{n} \boldsymbol{\Omega}^{\log} + \boldsymbol{\Omega}^{\log}\mathbf{S}_{\text{solid}}^{n}\right).
\end{equation}

\subsubsection{Plastic Correction}
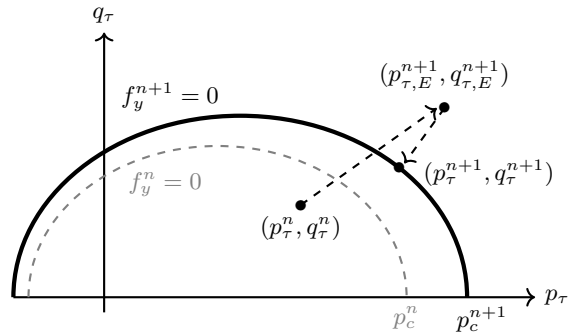
\begin{figure}[htb!]
    \centering 
    \begin{tikzpicture}[scale=1]
        \draw[thick, ->] (-0.2,0) -- (6.7,0) node[right] {\small $p_{\tau}$};
        \draw[thick, ->] (1,-0.2) -- (1,3.5) node[above] {\small $q_{\tau}$};
        
        \draw[thick, gray, dashed] (0,0) arc[start angle=180, end angle=0, x radius=2.5cm,
            y radius =2cm];
        \draw[ultra thick] (-0.2,0) arc[start angle=180, end angle=0, x radius=3cm,
            y radius =2.4cm];

        \node[right, gray] at (1.2,1.5) {\small $f_y^n = 0$};
        \node[right] at (1.1,2.65) {\small $f_y^{n+1} = 0$};

        \draw (3.6,1.2) node {\textbullet} node[below] {\small $(p_\tau^n,q_\tau^n)$};
        \draw[dashed, thick, ->] (3.6,1.22) -- (5.4,2.5);
        


        \draw (5.5,2.5) node {\textbullet};
        \draw[dashed, thick, ->] (5.5,2.5) -- (4.97,1.77);
        \node[above] at (5.5,2.6) {\small $(p_{\tau,E}^{n+1},q_{\tau,E}^{n+1})$};
        
        \draw (4.9,1.7) node {\textbullet};
        \node[right] at (5.1, 1.65) {\small $(p_\tau^{n+1},q_\tau^{n+1})$};
        
        \node[below, gray] at (5,0) {\small $p_c^n$};
        \node[below] at (6,0) {\small $p_c^{n+1}$};

    \end{tikzpicture}
    \caption{Illustration of elastic trialstep and reprojection. $(p_\tau^n,q_\tau^n)$ is the invariants pair of the old time step and $(p_{\tau,E}^{n+1},q_{\tau,E}^{n+1})$ the elastic trial step. $(p_\tau^{n+1},q_\tau^{n+1})$ is the final stress state, which was reprojected via the correction from equation~\eqref{eq:NewtonIteration} to the updated yieldsurface with the Newton iteration.}
    \label{fig:trialstep_and_reprojection}
\end{figure}

To evaluate the material state at each collocation point, the Cauchy pressure $p_E^{n+1}=\tilde{p}_{\text{dyn,E}}^{n+1}+p_{\text{hyd,E}}^{n+1}+c^{n+1}$ and the deviatoric stress $\mathbf{S}_{\text{solid,E}}^{n+1}$ are transformed into their respective Kirchhoff stress invariants, $p_{\tau,E}^{n+1}$ and $q_{\tau,E}^{n+1}$. If the trial state violates the yield condition $f_y(p_{\tau,E}^{n+1}, q_{\tau,E}^{n+1}) > 0$, we perform a plastic correction. In \cite{Blatny2024} it is shown that by expressing the updated stress invariants $p_{\tau}^{n+1}$ and $q_{\tau}^{n+1}$ as explicit functions of the plastic multiplier $\Delta \gamma$ (cf. equation~\eqref{eq:equivplasticstrain} ), the tensor return mapping reduces to finding the root of a single scalar residual equation
\begin{equation}
    R(\Delta \gamma) = f_y\left(p_{\tau}^{n+1}(\Delta \gamma), p_{\tau}^{n+1}(\Delta \gamma),p_{\tau,c}^{n+1}(\Delta\gamma) \right) = 0. \label{eq:NewtonIteration}
\end{equation}
This scalar equation is solved efficiently via a local Newton-Raphson iteration. An illustration is shown in figure~\ref{fig:trialstep_and_reprojection}.

Upon convergence of $\Delta \gamma$, we extract the coupling terms for the global solver 
\begin{equation}
    \mathbf{S}^{n+1}_{\text{solid}} =  \frac{q_\tau^{n+1}}{q_{\tau,E}^{n+1}}\mathbf{S}_{\text{solid,E}}^{n+1}, \quad\quad (\overline{\nabla\cdot \mathbf{v}})_P^{n+1} = -\frac{\Delta \gamma}{\Delta t}\frac{\partial f_y /\partial p}{||\partial f_y / \partial \bsigma||}
\end{equation}
(cf. equation~\eqref{eq:plasticflowrule_vol}) and update the internal hardening variable (cf. equation~\eqref{eq:pc_update})
\begin{equation}
    p_c^{n+1} = p_c^{n} \exp\left(-\xi\Delta t\left(\nabla\cdot \mathbf{v}\right)_P^{n+1} \right).
\end{equation}
Furthermore, the return mapping procedure yields an effective elastoplastic compressibility $\kappa_{\text{eff}}$. To ensure compliance with the conservation equations, we perform a post-processing step to globally resolve for $p_{\text{hyd}}$ and $p_{\text{dyn}}$, instead of propagating the Cam-Clay-corrected pressure into the subsequent time step. The updated compressibility, alongside the divergence of the deviatoric stress  $\nabla \cdot \mathbf{S}^{n+1}_{\text{solid}}$ is utilized to compute the final hydrostatic pressure
\begin{align}
-\frac{\kappa_{eff}}{\Delta t^2} p_{\text{hyd}}^{n+1}+\hat{\nabla}\cdot \left(\frac{{\rm 1}}{\rho } \hat{\nabla} p_{hyd}^{n+1} \right)=-\frac{\kappa_{eff}}{\Delta t^2} p_{\text{hyd}}^{n}+ \hat{\nabla}\cdot \left(\frac{{\rm 1}}{\rho } \left(\hat{\nabla}\cdot \mathbf{S}_{\text{solid}}^{n+1} \right)_{}^{T} \right)+ \hat{\nabla}\cdot \mathbf{g}.\label{eq:hydpressure2} 
\end{align}
Similar we obtain the final Poisson-type equation for the dynamic pressure by utilizing $(\overline{\nabla\cdot \mathbf{v}})^{n+1}=(\overline{\nabla\cdot \mathbf{v}})_P^{n+1}+(\overline{\nabla\cdot \mathbf{v}})_E^{n+1}$ as a source term
\begin{align}
-\frac{\kappa_{eff}}{\Delta t^2}p_{\text{dyn}}^{n+1}+\hat{\nabla}\cdot \left(\frac{{\rm 1}}{\rho } \hat{\nabla} p_{dyn}^{n+1} \right)=&-\frac{\kappa_{eff}}{\Delta t^2} p_{\text{dyn}}^{n}+\frac{1}{\Delta t}\left(\nabla\cdot v\right)^n-\frac{1}{\Delta t}\left(\overline{\nabla\cdot\mathbf{v}}\right)^{n+1}\nonumber\\&-\Phi \left(\mathbf{v}^{n+1}\right)+\hat{\nabla}\cdot \left(\frac{{\rm 1}}{\rho } \left(\hat{\nabla}\cdot \mathbf{S}_{\text{visc}}^{} \right)_{}^{T} \right).\label{eq:dynpressure2} 
\end{align}
Finally, the point positions are updated according to the velocity. This step incorporates a dedicated volume correction that explicitly accounts for the plastic divergence, thereby ensuring mass conservation during compaction. For a comprehensive discussion on general mass-conserving correction schemes within meshfree frameworks, the reader is referred to \cite{Suchde2023}.

\begin{figure}[htbp]
    \centering 
\begin{tikzpicture}[
    auto,
    block/.style={
        rectangle, 
        draw=black, 
        thick, 
        fill=white,
        minimum width=8.5cm,
        minimum height=1.5cm,
        inner sep=6pt,
        align=center, 
        font=\footnotesize\sffamily
    },
    resultBlock/.style={
        rectangle,
        rounded corners=4pt, 
        draw=black,
        thick,         
        dashed,        
        fill=black!10, 
        minimum width=5cm,
        minimum height=0.6cm,
        inner sep=4pt,
        align=center,
        font=\footnotesize\sffamily 
    },
    line/.style={
        draw, 
        -{Latex[length=3mm, width=2mm]}, 
        thick
    },
    shortLine/.style={
        draw,
        -{Latex[length=2mm, width=1.5mm]},
        thick 
    }
]
    \node [block] (corrTop) {
        \textbf{Global PDE (Elastic Pressure Trial)}\\[0.3em]
        {$\footnotesize\displaystyle \tilde{p}_E^{n+1} = \text{Poisson}\left(\mathbf{v}^{n},  \nabla\mathbf{S}_{\text{solid}}^{n+1}\right)$}
    };
    \node [resultBlock, below=0.15cm of corrTop] (corrTopRes) {
        \textbf{\textit{Solution:}}\hspace{0.2cm} $\footnotesize\displaystyle
        \begin{aligned}
            &\tilde{p}_E^{n+1}
        \end{aligned}$
    };
    \draw [shortLine] (corrTop) -- (corrTopRes);

    \node [block, below=0.5cm of corrTopRes] (pred) {
        \textbf{Global PDE (Elastic Trial)}\\[0.3em] 
        {$\footnotesize\displaystyle
\begin{aligned}
  \begin{lcases}
        \frac{\tilde{\mathbf{v}}^{n+1}-\mathbf{v}^{n}}{\Delta t}+ \frac{{\rm 1}}{\rho} \hat{\nabla} c^n+1-\frac{{\rm 1}}{\rho } \left(\hat{\nabla}\cdot \mathbf{S}_{\text{visc}}\left(\mathbf{v}^{n+1}, \hat{\eta}\right) \right)_{}^{T} = &-\frac{{\rm 1}}{\rho } \hat{\nabla} \tilde{p}_E^{n+1}+\frac{{\rm 1}}{\rho } \left(\hat{\nabla}\cdot \mathbf{S}_{\text{solid}}^{n} \right)_{}^{T}+\mathbf{g}\\
        -\frac{\kappa_E}{\Delta t}c^{n+1}+ \hat{\nabla}\cdot\left(\frac{\Delta t_{\text{virt}}}{\rho} \hat{\nabla} c^{n+1}\right)- \hat{\nabla}\cdot\tilde{\mathbf{v}}^{n+1} =&-(\overline{\nabla\cdot\mathbf{v}})_E&
    \end{lcases}
\end{aligned}$ }
    };
    \node [resultBlock, below=0.15cm of pred] (predRes) {
        \textbf{\textit{Solution:}}\hspace{0.2cm} $\footnotesize\displaystyle
        \begin{aligned}
            \mathbf{\tilde{\tilde{v}}}^{n+1} &= \mathbf{\tilde{v}}^{n+1} -\frac{\Delta t_{\text{virt}}}{\rho }\nabla c, \hspace{0.4cm} p_E^{n+1} &=\tilde{p}_E^{n+1} + c
        \end{aligned}$ 
    };
    \draw [line] (corrTopRes) -- (pred);
    \draw [shortLine] (pred) -- (predRes);

    \node [block, below=0.5cm of predRes] (point) {
        \textbf{Pointwise Operation (Elastic Trial)}\\[0.3em]
        {$\footnotesize\displaystyle \begin{aligned}
    \mathbf{S}_{\text{solid,E}}^{n+1} =\mathbf{S}_{\text{solid}}^{n} + \Delta t \left(2 \mu_E \text{dev}\left(\mathbf{D}^{n+1}\right) - \mathbf{S}_{\text{solid}}^{n} \boldsymbol{\Omega}^{\log} + \boldsymbol{\Omega}^{\log}\mathbf{S}_{\text{solid}}^{n}\right) \end{aligned}$}
    };
    \node [resultBlock, below=0.15cm of point] (pointRes) {
        \textbf{\textit{Solution:}}\hspace{0.2cm} $\footnotesize\displaystyle
        \begin{aligned}
            \mathbf{S}_{\text{solid,E}}^{n+1} 
        \end{aligned}$
    };
    \draw [line] (predRes) -- (point);
    \draw [shortLine] (point) -- (pointRes);

    \node [block, below=0.5cm of pointRes] (CCC) {
        \textbf{Pointwise Operation (MCC Plastic Correction)}\\[0.3em]
        {$\footnotesize\displaystyle\begin{aligned} 
    \begin{cases}
        \text{Reprojection} & \text{if } f_y(p_{E}^{n+1},\mathbf{S}_{\text{solid},E}^{n+1}, p_c^{n}) > 0 \\
        \text{Keep elastic Solution} & \text{else}
    \end{cases}\end{aligned}$}
    };
    \node [resultBlock, below=0.15cm of CCC] (CCCRes) {
       \textbf{\textit{Solution:}}\hspace{0.2cm} $\footnotesize\displaystyle
        \begin{aligned}
            \nabla\mathbf{S}_{\text{solid}}^{n+1}, \hspace{0.2cm} \left(\overline{\nabla\cdot\mathbf{v}} \right)_{P}^{n+1}, \hspace{0.2cm} p_c^{n+1}
        \end{aligned}$
    };
    \draw [line] (pointRes) -- (CCC);
    \draw [shortLine] (CCC) -- (CCCRes);

    \node [block, below=0.5cm of CCCRes] (corr) {
        \textbf{Global PDE (Pressure-Corrector)}\\[0.3em]
        {$\footnotesize\displaystyle p^{n+1} = \text{Poisson}\left(\tilde{\tilde{\mathbf{v}}}^{n+1},  \nabla\mathbf{S}_{\text{solid}}^{n+1},  \left( \overline{\nabla\cdot\mathbf{v}} \right)_{P}^{n+1}\right)$}
    };
    \node [resultBlock, below=0.15cm of corr] (corrRes) {
        \textbf{\textit{Solution:}}\hspace{0.2cm} $\footnotesize\displaystyle
        \begin{aligned}
            &p^{n+1}, \hspace{0.2cm} \mathbf{v}^{n+1} \hspace{0.1cm} \left(\text{Correction based on } \left( \overline{\nabla\cdot\mathbf{v}} \right)_{P}^{n+1} \right)
        \end{aligned}$
    };
    \draw [line] (CCCRes) -- (corr);
    \draw [shortLine] (corr) -- (corrRes);

    
    \draw [line] ($(corrTop.north)+(0,0.8)$) -- (corrTop.north) node[midway, right, font=\small\bfseries] {Start ($t_0$)};

    \path (pred.west) ++(-0.8, 0) coordinate (safeLeft);

    \draw [line] (corrRes.south) 
        -- ++(0,-0.5) coordinate(bottomPoint) 
        -- node[below, pos=0.5, font=\small\bfseries] {Next Timestep ($t_{n+1}$)} 
           (bottomPoint -| safeLeft) 
        |- (corrTop.west);

\end{tikzpicture}
    \caption{Scheme of the numerical solver.}
    \label{fig:Solver-Scheme}
\end{figure}
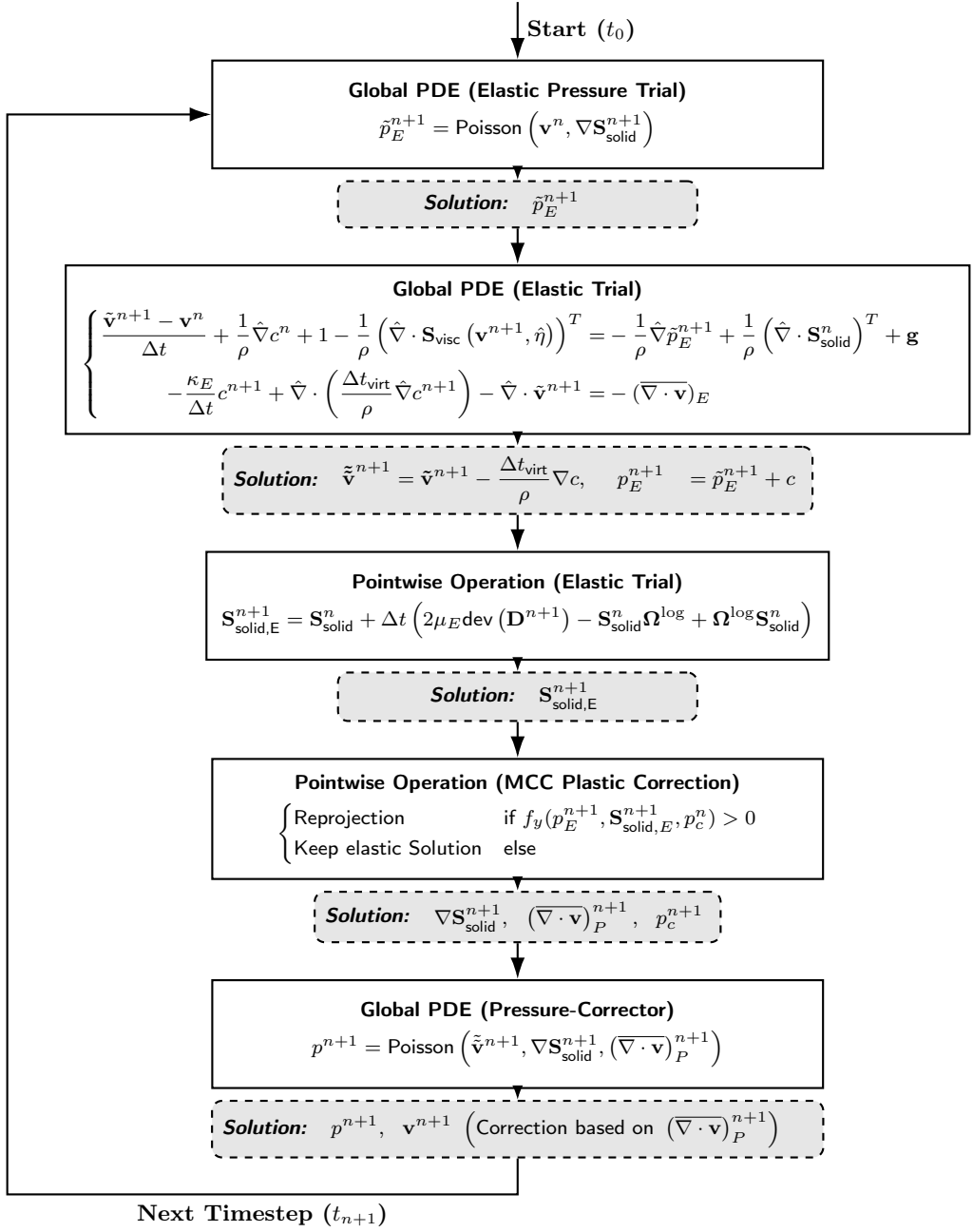

An overview of the scheme is shown in figure~\ref{fig:Solver-Scheme}.

\subsection{Advantages of Snow in GFDM}
The \ac{gfdm} offers several general advantages over numerical methods previously utilized for snow simulations. As a strong-form approach, it avoids the strict polynomial order requirements typically associated with weak-form methods, such as the Ladyzhenskaya-Babuška-Brezzi conditions.

In this framework, a computational node does not carry mass (see, for example, \cite{Suchde2023}). Instead, it functions exclusively as a collocation point, where the governing partial differential equations are evaluated. A primary benefit of such collocation point-based meshfree methods is the simplicity of modifying the point cloud. To address domain distortion caused by material movement, \ac{gfdm} points can be seamlessly added to fill artificial holes, or merged and deleted when they become too closely packed \cite{Suchde2019, suchde2023point}. This flexibility facilitates straightforward, continuous adaptive refinement. In contrast, while multi-resolution frameworks for particle-based methods are being developed, they generally remain limited to a few discrete levels of refinement, whereas collocation methods have long supported continuous resolution adaptivity.

Another crucial advantage is discrete consistency: the procedure for calculating discrete derivatives in \ac{gfdm} guarantees consistency (usually ensuring polynomial reproducibility up to a targeted order) at both continuous and discrete levels \cite{halada2025overview}. In particle-based methods, achieving this discrete consistency is challenging because it heavily relies on the specific particle arrangement and is difficult to maintain for arbitrary distributions. Additionally, the method operates entirely without a background grid, eliminating the need to project data between grids and particles, and it generally permits larger stable time steps than \ac{mpm}, comparable to those of \ac{sph}.

Finally, meshfree collocation approaches allow for the easy enforcement of a wide variety of boundary conditions, a process that is typically far more complex in particle-based methods like \ac{sph} and \ac{mpm}.  

A more detailed discussion of these methodological advantages can be found, for example, in \cite{Suchde2023}.

\subsection{Boundary Conditions}

A challenge in Lagrangian particle methods is the accurate imposition of boundary conditions. Formulations like \ac{sph} suffer from so-called particle deficiency (kernel truncation) at boundaries, degrading consistency and requiring artificial strategies like ghost particles. Similarly, the \ac{mpm} can introduce grid-dependency artifacts when resolving complex curved geometries, which limits the contact precision.

The \ac{gfdm} overcomes these limitations by operating as a strong formulation collocation scheme. It does not rely on integral approximations but reconstructs differential operators directly. This allows boundary points to be treated as collocation nodes, enabling the direct enforcement of Dirichlet and Neumann conditions without approximation deficits.

For snow mechanics, this framework offers advantages in interaction with complex curved boundaries (necessary for vehicle applications), which are represented exactly by the boundary nodes and avoiding the grid alignment errors like in \ac{mpm}. On top of that, as the stress tensor is evaluated at the boundary, application of direct physical friction laws (e.g., Coulomb friction) based on the pressure is possible, instead of relying on numerical penalty forces common in other methods.

Since boundary conditions depend on the application and strongly influence material behavior, we will discuss their explicit formulation in the corresponding chapter of the simulation results.

    \section{Simulation Results}\label{sec:results}

We apply our model to different benchmark problems, and validate the model against standard benchmarks for snow and free surface flow problems. We also show convergence of our simulations.
First, the hardening and compaction behavior is validated against a hydrostatic cyclic compression test from the literature, accompanied by a rigorous spatial and temporal convergence analysis. Next, we benchmark the model against a standard shear box test, addressing the well-known strain localization phenomenon associated with material softening. Finally, we demonstrate the \ac{fsi} capabilities through a sphere drop test into a snow-filled box. To underscore the broader potential of the proposed framework, the section concludes with a demonstrative, qualitative \ac{fsi} simulation of a full vehicle interacting with snow. All simulations were performed using the software \cite{MESHFREE}, in which the described modeling and scheme are implemented.

\subsection{Cyclic Hydrostatic Compression}

We begin with a hydrostatic compression test to validate the theoretical framework, presented in section~\ref{subsec:cam-clay}, against the experimental data from \cite{Meschke1996}. For the 3D simulation, a cubic domain is filled with collocation points, and slip boundary conditions are applied at the boundaries. Specifically, these conditions are governed by 
\begin{equation} 
\mathbf{n}\cdot \nabla p_{\text{hyd}} = \rho \left( \mathbf{n} \cdot \mathbf{g}   \right) + {\bf n} \cdot \left( \nabla \cdot{\bf S}_{\text{solid}} \right)
\end{equation}
for the hydrostatic pressure where $\mathbf{n}$ is the boundary normal of the geometry walls. It should be noted that the hydrostatic pressure component was deactivated for this particular experiment. In a closed domain characterized by complex boundaries and purely Neumann conditions, the kinematics of the collocation points are governed entirely by the pressure gradients rather than absolute pressure values. Therefore, the dynamic pressure represents the total pressure, avoiding the ill-posed nature of a pure Neumann problem for the Poisson equation. For the dynamic pressure the slip boundary condition reads
\begin{equation} \mathbf{n}\cdot\nabla p_\text{dyn}= -\rho{\bf n} \cdot \left(   \left( ( {\bf v}-{\bf v}_\text{wall}) \cdot \nabla \right) ({\bf v}-{\bf v}_\text{wall}) +  \frac{D}{Dt} {\bf v}_\text{wall}   \right) + {\bf n} \cdot \left( \nabla\cdot {\bf S}_{\text{visc}} \right),
\end{equation}
where $\mathbf{v}_{wall}$ denotes the velocity of the geometry walls. The equations can be simply derived by multiplying the boundary normal from the left \cite{MESHFREE_BND_wall}. The velocity of the boundary points is given by a Dirichlet condition ($\mathbf{v}\cdot \mathbf{n}=\mathbf{v}_{wall}\cdot \mathbf{n}$) in normal direction and a Neumann condition in tangential direction of the wall. Consequently, frictional effects, which are neglected under the current slip boundary conditions, can be incorporated into the model via the deviatoric solid stress term $ {\bf S}_\text{solid}$. For the simplified case under consideration, the normal deviatoric stress components are set to zero in tangential directions:
\begin{align} 
  {\bf t}_1 \cdot \bf \Ss \cdot {\bf n} = 0 \\ 
  {\bf t}_2 \cdot \bf \Ss \cdot {\bf n} = 0
\end{align}
where $\bf t_1$ and $\bf t_2$ are the tangentials associated to $\bf n$, i.e. 
\begin{align} 
  {\bf t}_1 \cdot {\bf n} = 0 \\ 
  {\bf t}_2 \cdot {\bf n} = 0
\end{align}
such that
\begin{align} 
  \left\Vert{\bf{t}_1}\right\Vert = \left\Vert{\bf{t}_2}\right\Vert = \left\Vert{\bf n}\right\Vert =1.
\end{align}
For the sake of brevity, the explicit mathematical formulations of the boundary conditions are omitted in the subsequent sections. Conceptually, no-slip and frictional constraints can be imposed analogously via the deviatoric stress tensor. At free surfaces, the roles of the primary variables are reversed: a Dirichlet boundary condition is prescribed for the pressure, while a Neumann boundary condition is applied to the velocity field. For a comprehensive description of the boundary treatment within this framework, the reader is referred to \cite{MESHFREE_BND_wall}.

In this setup, three walls are kept stationary, while the remaining three walls move cyclically at a uniform velocity toward the centroid of the cube. Compare with the following Figure~\ref{fig:cube-kobs} and Figure~\ref{fig:cyclic-compress-experiment}. The temporal evolution of the density of the snow cube, alongside the moving wall geometries, is further illustrated in Supplementary Video 1.
\begin{figure}[htbp]
    \centering
    \includegraphics[width=1.0\textwidth]{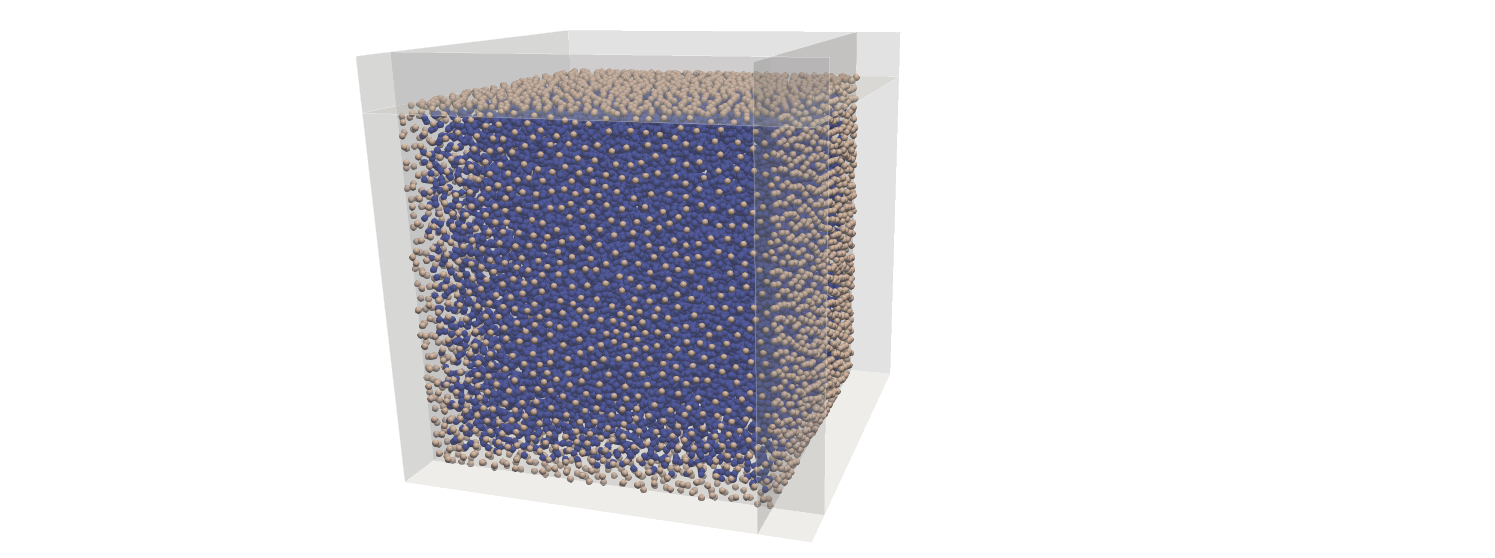} 
    \caption{Classification of collocation points in the cyclic compression test: Inner points (blue) and boundary points (off-white) of the transparent geometry. Three faces of the cubic domain are subjected to a cyclic inward and outward displacement.}
    \label{fig:cube-kobs}
\end{figure}
\begin{figure}[htbp]
    \centering
    \includegraphics[width=0.8\textwidth]{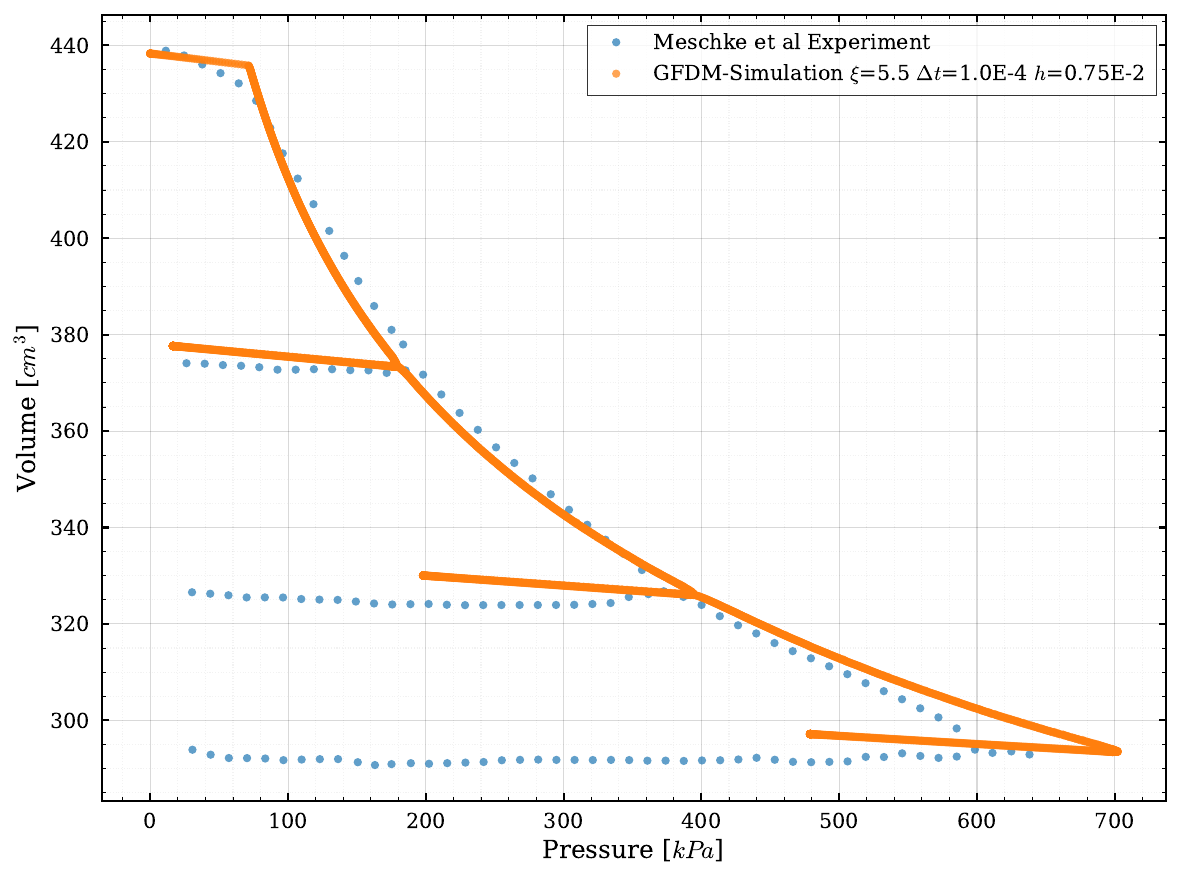} 
    \caption{Comparison of the GFDM simulation with the experimental hydrostatic cyclic compression test from \cite{Meschke1996}.}
    \label{fig:cyclic-compress-experiment}
\end{figure}

\begin{figure}[htbp]
    \centering
    \includegraphics[width=0.8\textwidth]{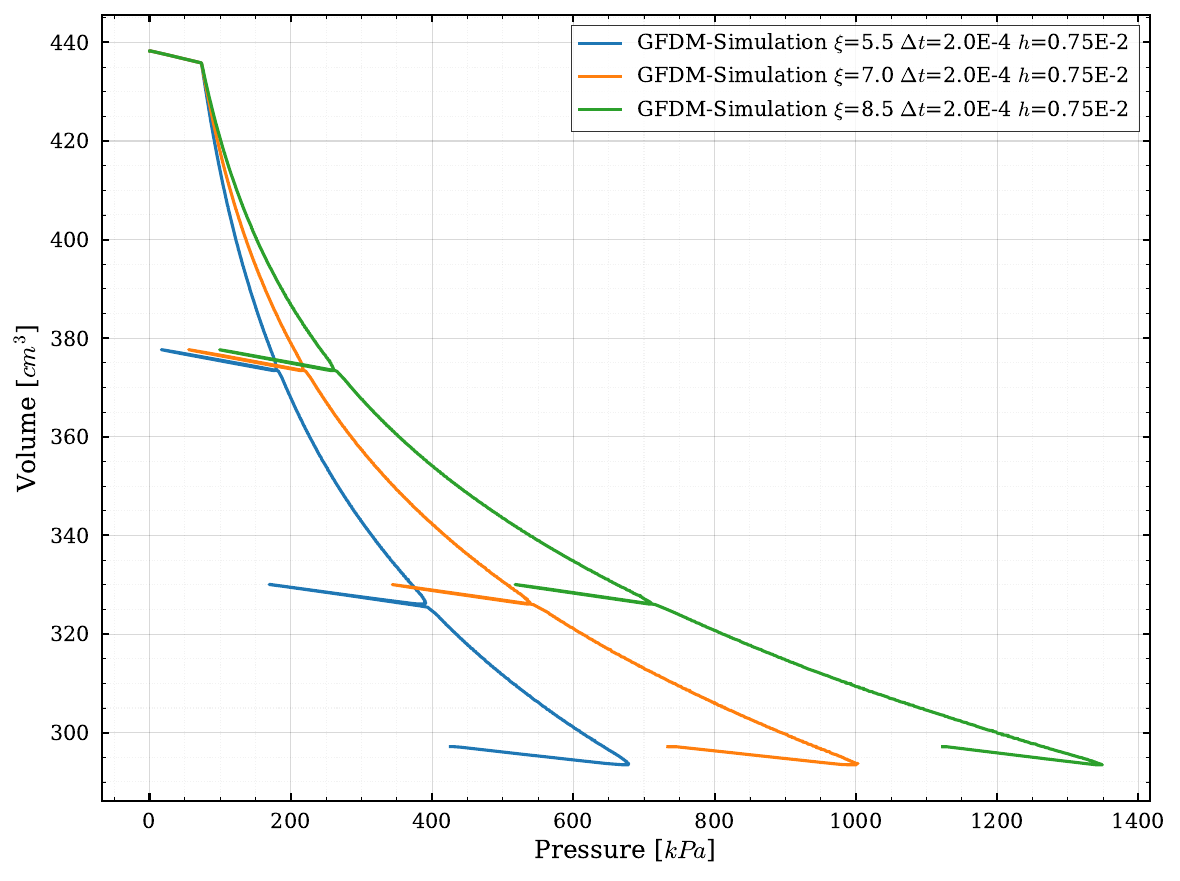} 
    \caption{Influence of the hardening parameter $\xi$ in the GFDM simulation for the cyclic compression test.}
    \label{fig:cyclic-compress-xi}
\end{figure}

Furthermore, the subsequent Figure~\ref{fig:cyclic-compress-xi} illustrates the significant impact of the hardening parameter on the mechanical behavior of snow. As one would expect, increasing the hardening parameter $\xi$ yields a steeper slope of the virgin compression line. Consequently, under identical loading conditions, higher pressure levels are required to achieve the same degree of compression.

\begin{figure}[htbp]
    \centering
    \includegraphics[width=0.8\textwidth]{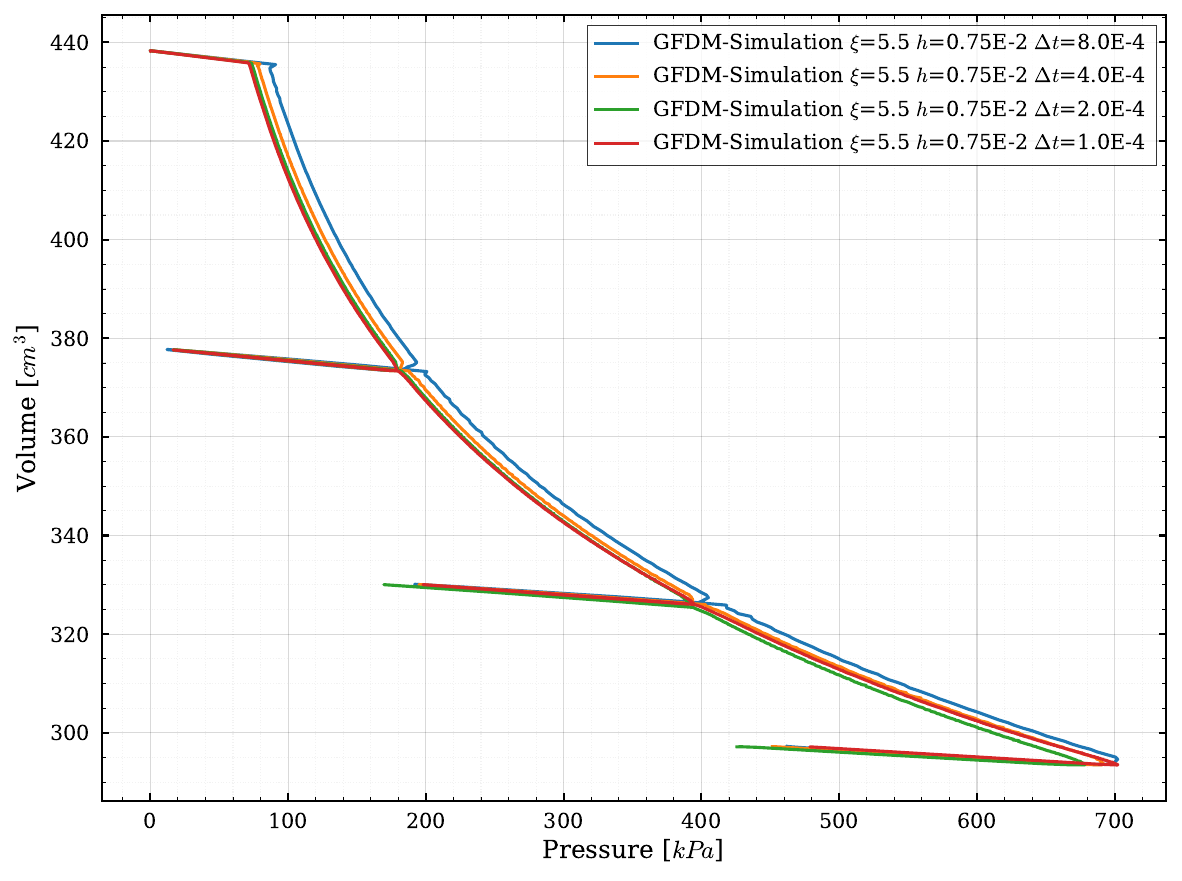} 
    \caption{Time step refinement of the GFDM simulation for the cyclic compression test.}
    \label{fig:cyclic-compress-refine_dt}
\end{figure}

\begin{figure}[htbp]
    \centering
    \includegraphics[width=0.8\textwidth]{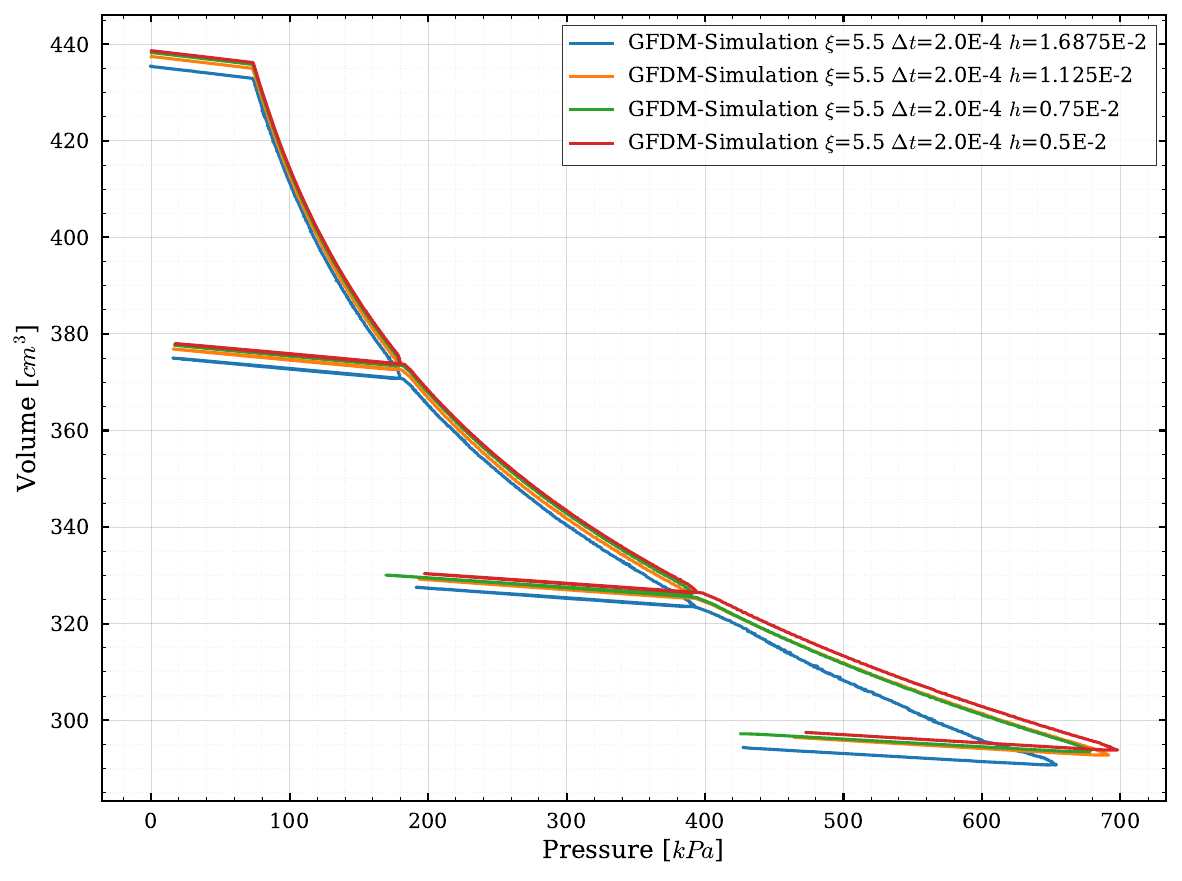} 
    \caption{Spatial refinement with respect to the smoothing length $h$ of the GFDM simulation for the cyclic compression test.}
    \label{fig:cyclic-compress-refine_h}
\end{figure}

For a specific hardening parameter $\xi$, time step size $\Delta t$, and smoothing length $h$, the results show agreement with the experiment. Furthermore, Figure~\ref{fig:cyclic-compress-refine_dt} and Figure~\ref{fig:cyclic-compress-refine_h} demonstrate convergence upon spatial and temporal refinement, respectively.

\begin{figure}[htbp]
    \centering
    \includegraphics[width=0.8\textwidth]{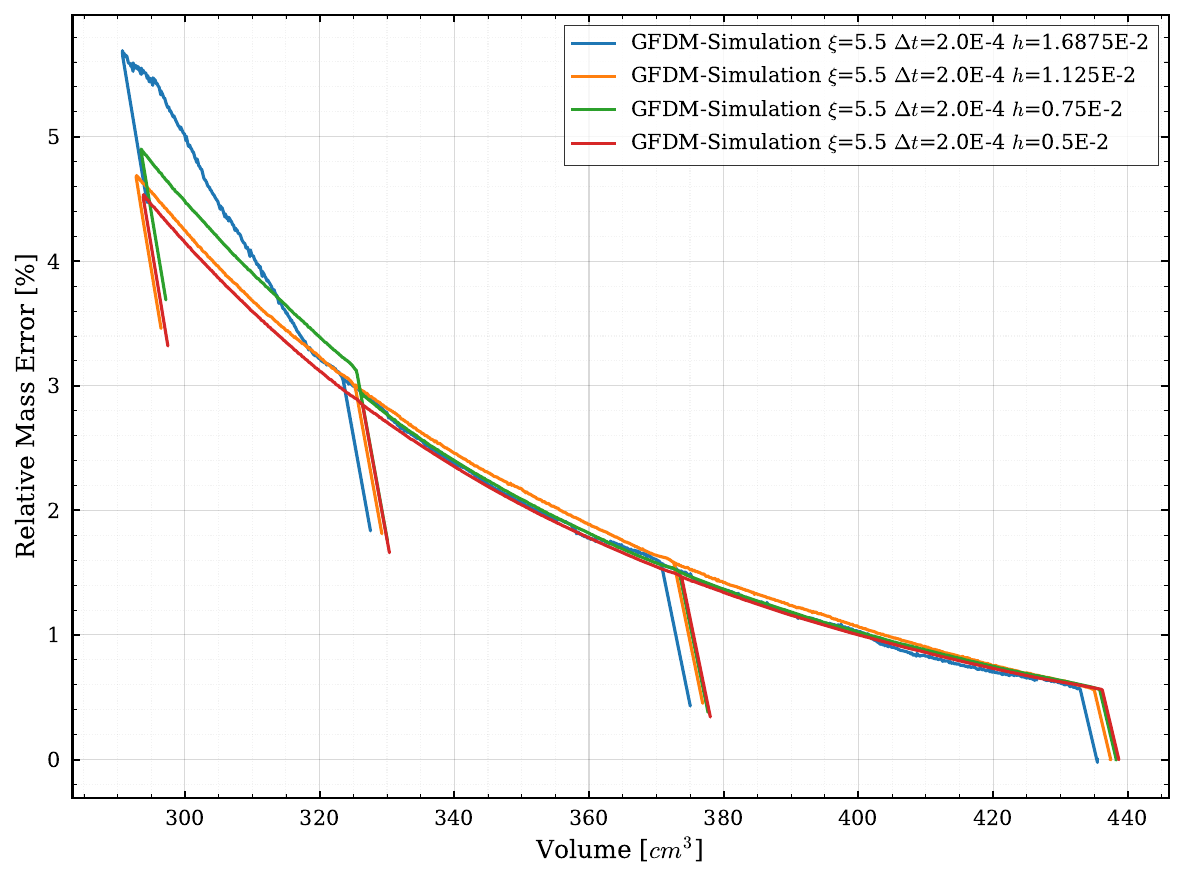} 
    \caption{Relative mass error, calculated from the initial mass and the current mass of the cube obtained via volume integration of the local density $\rho$.}
    \label{fig:cyclic-compress-mass}
\end{figure}

However, this convergence is strictly limited to the case of pure compression. As will be shown in the subsequent chapter, softening induces localization, which precludes convergence under mixed stress states. Finally, we demonstrate that by employing the correction method proposed in \cite{Suchde2023}, a volume reduction of $35\%$ yields a relative mass conservation error ($RME$) of $3\%$. This mass error is calculated according to 
\begin{equation}
    RME = \frac{\int_C \rho dV - m_{\text{init}}}{m_{\text{init}}}
\end{equation}
where $C$ it the snow cube domain and $m_\text{init}$ the initial mass of the snow cube. It can be observed in Figure~\ref{fig:cyclic-compress-mass} that spatial refinement leads to a corresponding reduction in the numerical error. This improvement is partially attributable to the enhanced accuracy of the volume integration. It should be noted that this represents an extreme case, where boundary effects, volume calculation derived from the point distribution, point deletion, and the absence of a free surface all exacerbate the issue of this error. Consequently, in free surface simulations, wherein correction algorithms exhibit enhanced efficiency, as well as in scenarios with reduced global relative volumetric variations, the associated relative mass error is rendered negligible. It is also important to note here that mass conservation errors occur across most Lagrangian discretization methods \cite{Suchde2023}.

\subsection{Shear Box Test}

\begin{figure}[htbp]
    \centering
    \includegraphics[width=1.0\textwidth]{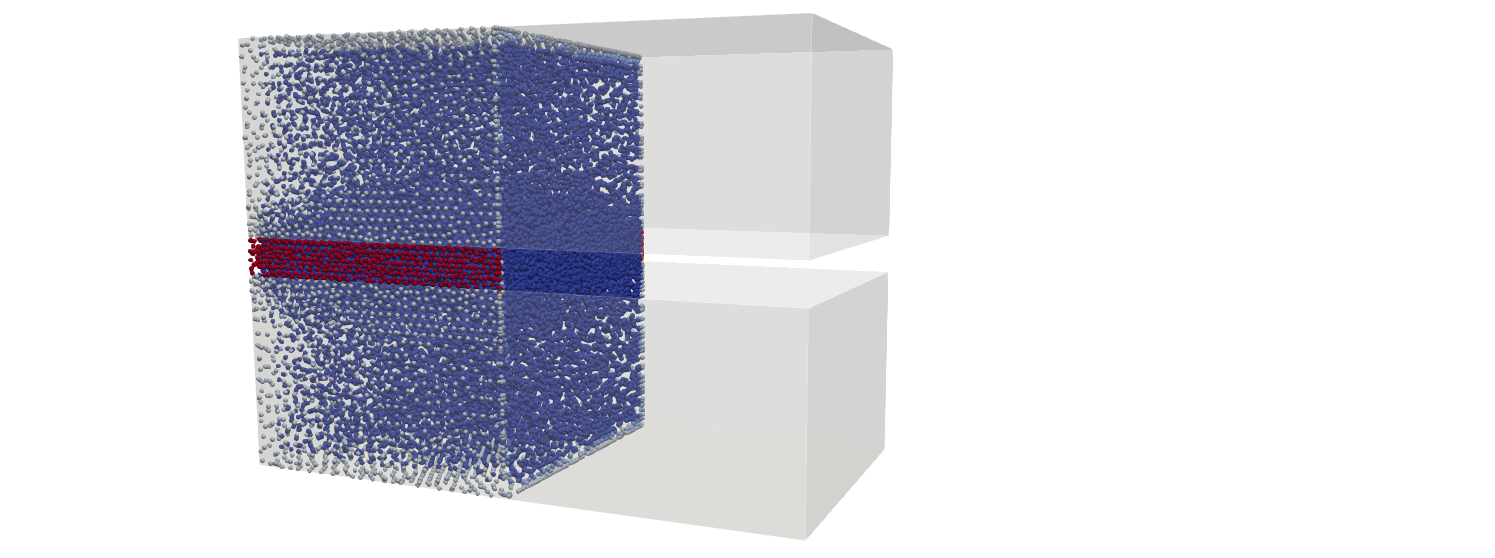} 
    \caption{Classification of collocation points in the shear box simulation: Inner points (blue), boundary points (grey) of the shear box geometry and free surface points (red) of the gap.}
    \label{fig:shearbox-kob}
\end{figure}

The focus now shifts to the shear box benchmark, adapted also from \cite{Meschke1996}. The objective is to demonstrate that while the softening model functions as intended, it encounters strain localization issues during spatial refinement. The computational domain is again discretized using collocation points (compare with Figure~\ref{fig:shearbox-kob}), deliberately leaving a $4 mm$ gap to serve as a free surface. No-slip boundary conditions are imposed on all remaining boundaries.

\begin{figure}[htbp]
    \centering
    \includegraphics[width=1.0\textwidth]{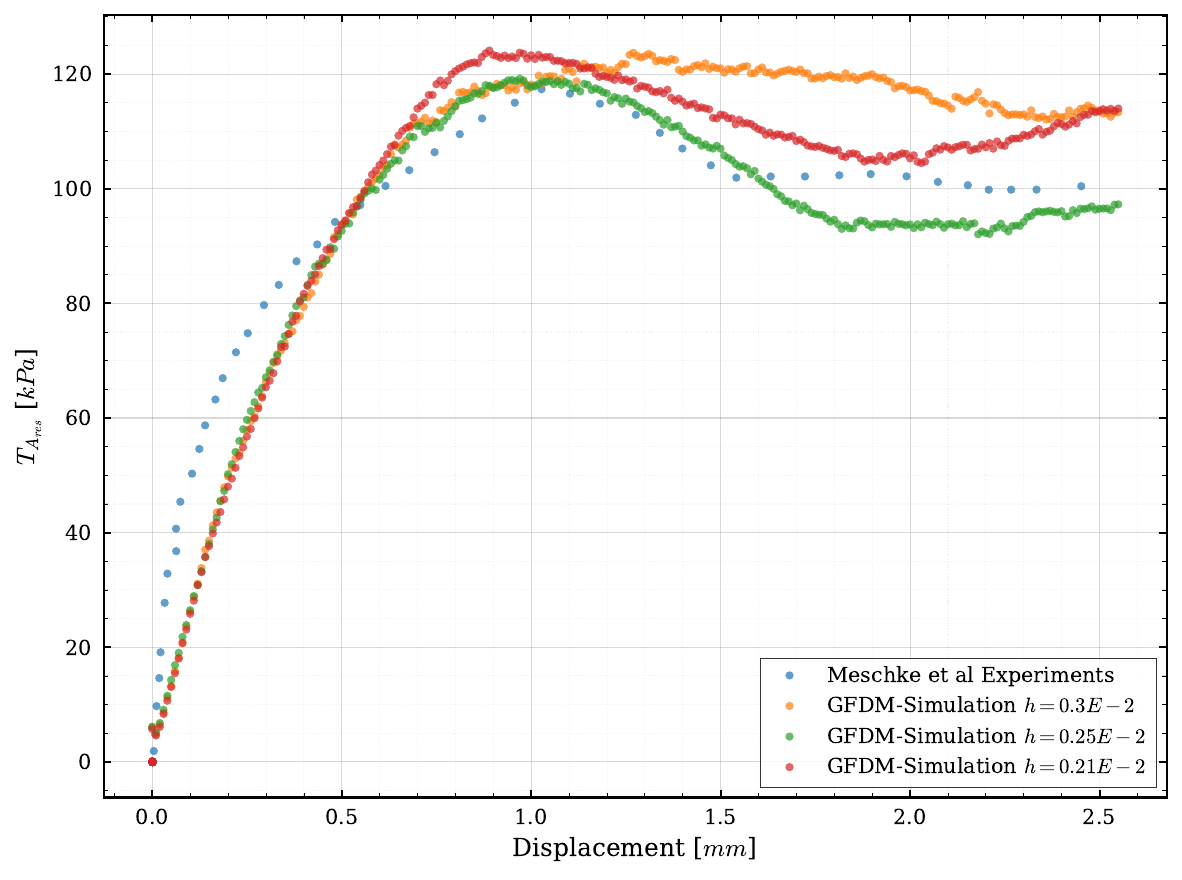} 
    \caption{Comparison of the GFDM simulation with the experimental shear box test from \cite{Meschke1996}.}
    \label{fig:shearbox-experiment}
\end{figure}

\begin{figure}[htbp]
    \centering
    \includegraphics[width=1.0\textwidth]{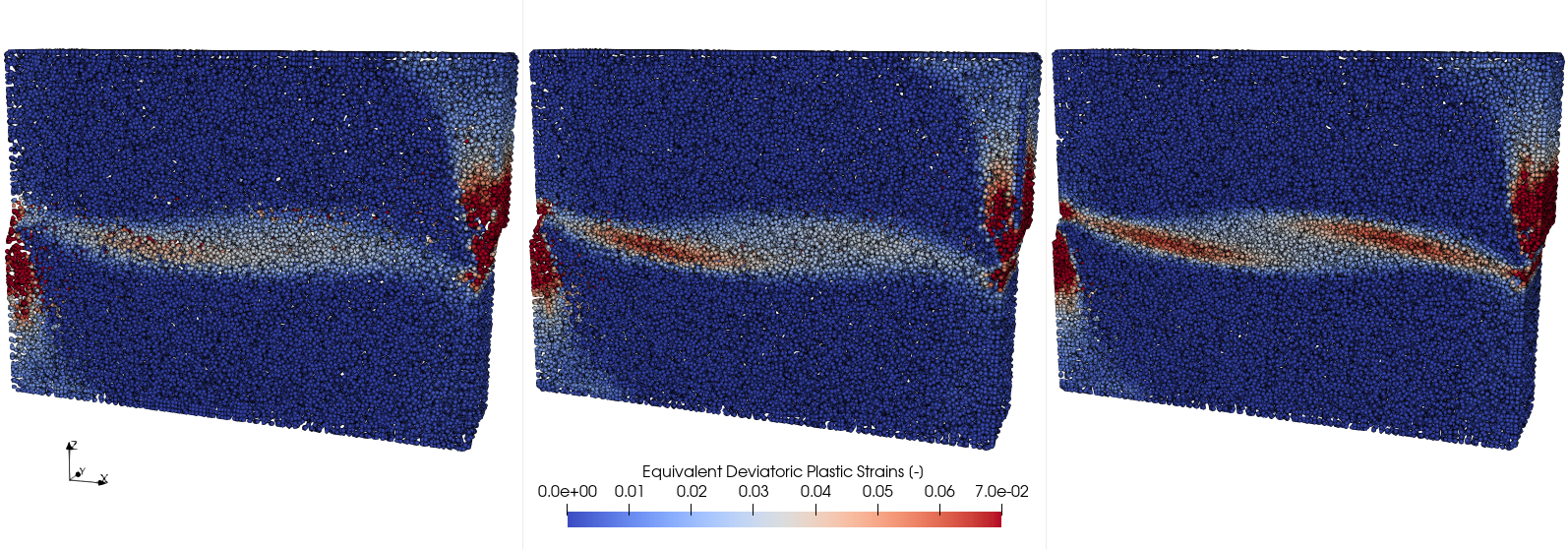} 
    \caption{Cross-sectional view of the shear band localization in the shear box test during spatial refinement of the gap for different smoothing lengths: $h=3.0\times 10^{-3}$ (left), $h=2.5\times 10^{-3}$ (middle), and $h=2.1\times 10^{-3}$ (right).}
    \label{fig:shearbox-localization}
\end{figure}

Interpreting the experimental documentation, this specific gap size is assumed to closely approximate the physical setup employed by \cite{Meschke1996}. It should be noted that maintaining a shear gap on the order of the characteristic grain size is standard convention in such setups. The macroscopic shear stress $T_{A_{res}}$ was subsequently determined by normalizing the horizontal shear force applied to the upper box by the transient residual contact area shared with the lower box. The experimental results obtained by \cite{Meschke1996} are plotted against our numerical predictions in Figure~\ref{fig:shearbox-experiment}.
The simulations contrast a uniform initial smoothing length $h=0.3\times10^{-2}$, with localized spatial refinements, $h=0.25\times10^{-2}$ and $h=0.21\times10^{-2}$, around the shear band. The load-displacement curves clearly demonstrate an artificial shift in the peak force and a significant alteration in the residual stress. Mechanically, this confirms the loss of structural objectivity inherent to local softening. As the GFDM resolution increases, the dissipated energy approaches zero, altering the macroscopic structural response. Figure~\ref{fig:shearbox-localization} visually corroborates this pathological strain localization. 
Supplementary Video 2 provides a dynamic visualization of the shear band evolution and the moving boundaries during the shear box test simulation. With increasing spatial resolution, the localization bands collapses into a more and more narrow crack, causing the deviatoric plastic strains to accumulate to locally extreme values. To overcome these numerical artifacts in future research, the model must be regularized and augmented with a suitable damage formulation (for example, \cite{Moeineddin2024, Moeineddin2026}). Crucially, this necessitates embedding a characteristic internal length scale, proportional to the physical grain size of the material, to regularize the evolution of plastic strains and ensure a resolution-independent fracture topology.

\subsection{Snow Structure Interaction (Sphere Drop)}\label{sec:spheredrop}
To evaluate the \ac{fsi} capabilities of the proposed framework, we investigate a classic general free-surface benchmark problem: the impact of a solid sphere dropped into a container filled with a granular medium. It should be noted that direct experimental validation data for this specific configuration with snow is currently unavailable in the literature. Nevertheless, the numerical results exhibit strong phenomenological similarities to the experimental sphere drop tests conducted with sand, as reported by \cite{Hu2021}. From a materials modeling perspective, this qualitative agreement in the macroscopic flow behavior demonstrates the physical plausibility of the model. Consequently, the primary objective of this test case is not strict experimental validation, but rather to illustrate the robustness and general feasibility of the FSI coupling within the present meshfree approach.

The ordinary differential equations governing the motion of the sphere's \ac{cog} are formulated as
\begin{align} 
\frac{d }{dt_{RB}} \left( \mathbf{x}_{COG} \right) &= \mathbf{v}_{COG} \nonumber\\
 \frac{d }{dt_{RB}} \left( m  \mathbf{v}_{COG} \right) &= \mathbf{F}_{snow} + \mathbf{F}_{gravity}  \\
 \frac{d}{dt_{RB}} \left( \mathbf{I} \boldsymbol{\omega}_{COG} \right) &= \mathbf{M}_{snow}, \nonumber
\end{align}
where $\mathbf{x}_{COG}$ are the coordinates of the \ac{cog}, $\mathbf{v}_{COG}$ the corresponding velocity, $\omega_{COG}$ the angular velocity and $\Delta t_{RB} \leq \Delta t$ the timestep size of the rigid body solver. $m$ and $\mathbf{I}$ are the mass respectively the rotational inertia tensor of the rigid body sphere. Here, $\mathbf{F}_{snow}+\mathbf{F}_{gravity}$ represents the resultant force exerted on the sphere, where $\mathbf{F}_{gravity}=m \mathbf{g}$ is the force due to gravity and $\mathbf{F}_{snow}$ is computed by integrating the stress tensor, evaluated at the boundary points of the snow domain, over the submerged surface area of the sphere. $\mathbf{M}_{snow}$ is the corresponding moment about the \ac{cog}.

\begin{figure}[htbp]
    \centering
    \includegraphics[width=0.8\textwidth]{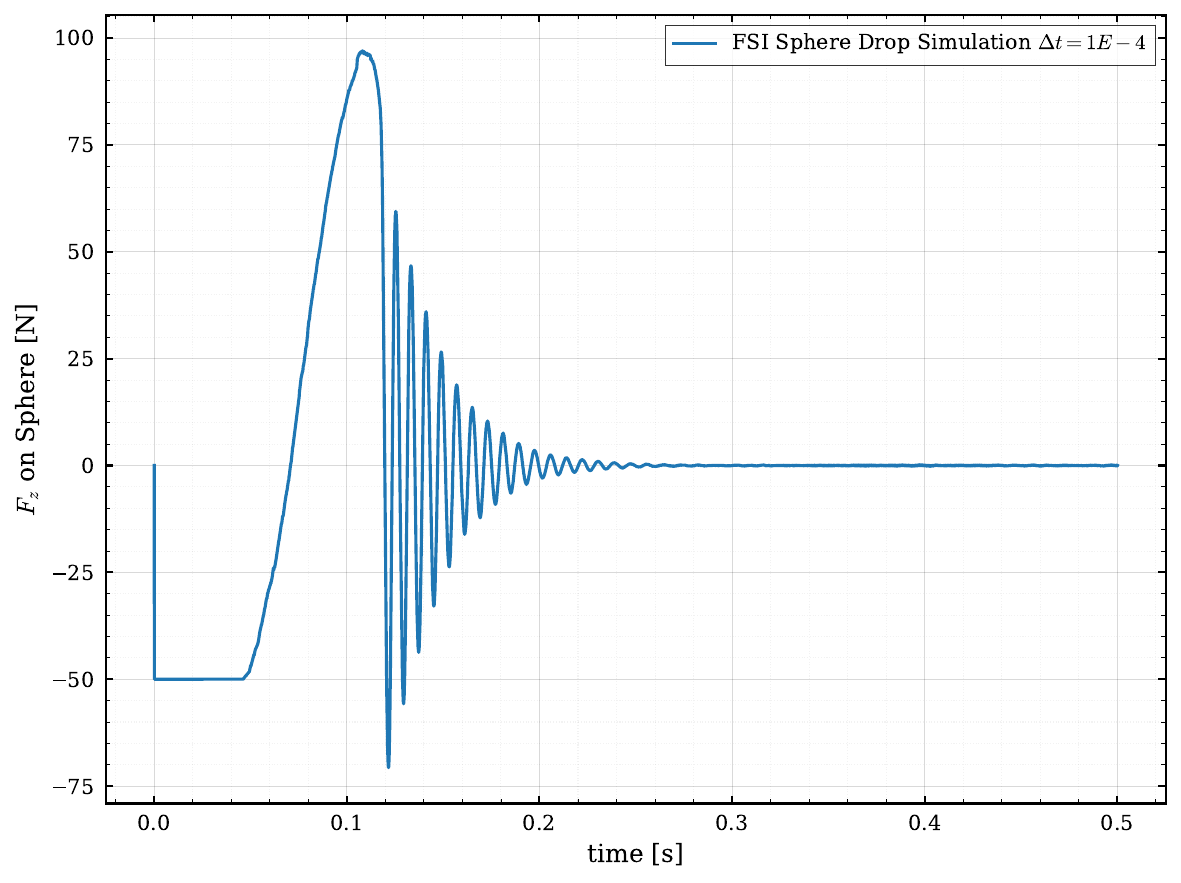} 
    \caption{Time history of the vertical force $\mathbf{F}_z=(\mathbf{F}_{fluid}+\mathbf{F}_{gravity})_z$ acting on the sphere.}
    \label{fig:spheredrop-Fz}
\end{figure}

Figure~\ref{fig:spheredrop-Fz} displays the temporal evolution of the vertical force $F_z=(\mathbf{F}_{fluid}+\mathbf{F}_{gravity})_z$ acting on the sphere. The rapidly decaying oscillations observed shortly before the body reaches its static equilibrium correspond to the elastic response of the snow model. It should be noted that, at present, no explicit physical viscosity model has been incorporated to provide further damping, additional to the inherent numerical viscosity, proportional to $\mu \Delta t$. Rather than indicating instability, these oscillations highlight the high resolving power of the numerical scheme and serve as a strong indicator of an accurately captured, purely elastoplastic material behavior.

\begin{figure}[htbp]
    \centering
    \includegraphics[width=0.8\textwidth]{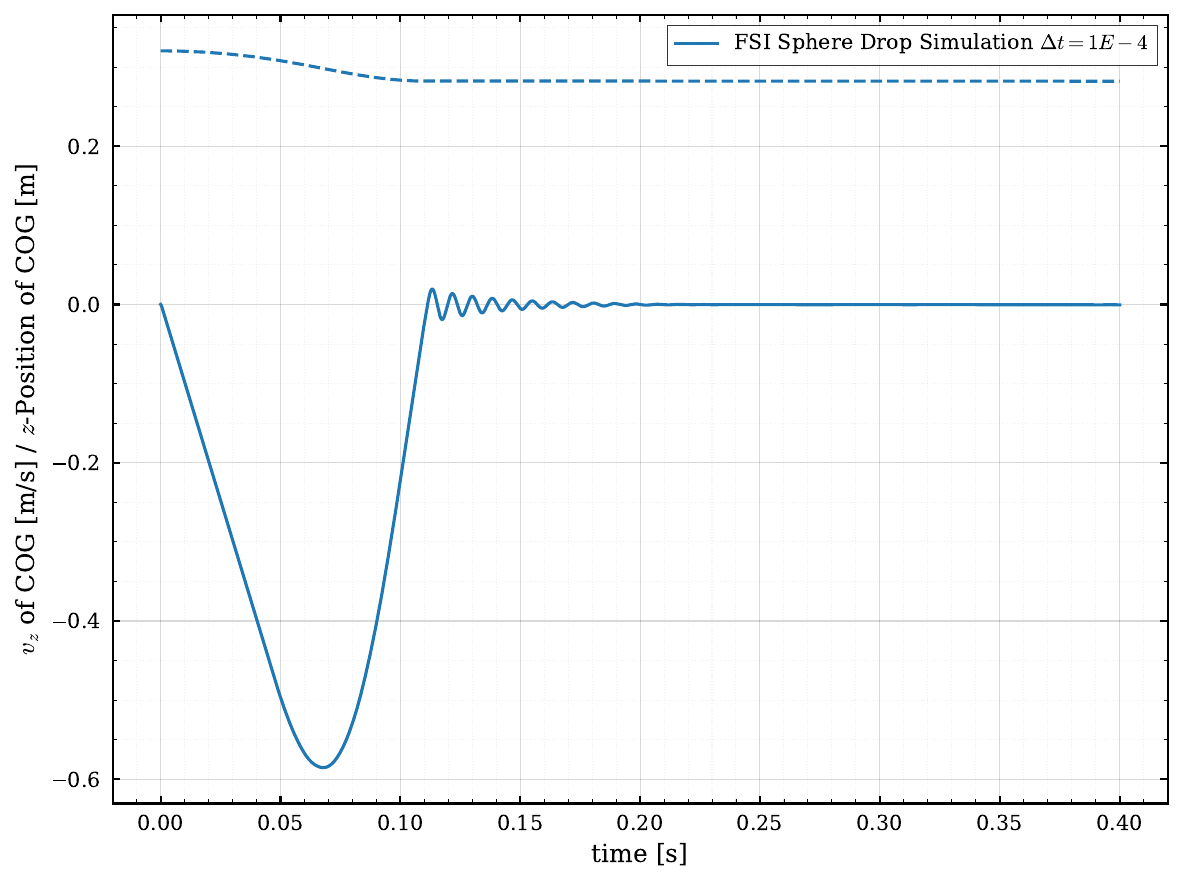} 
    \caption{Time history of the velocity ($v_z$) and the vertical position ($z$-direction) of the center of gravity of the sphere.}
    \label{fig:spheredrop-vz}
\end{figure}

Furthermore, the amplitude of the corresponding velocity oscillations in Figure~\ref{fig:spheredrop-vz} diminishes significantly as the sphere approaches its resting state, showing good qualitative agreement with the findings reported by \cite{Hu2021}. This stability is further corroborated by the vertical $z$-position of the sphere's \ac{cog}, also illsutrated in Figure~\ref{fig:spheredrop-vz}. The trajectory transitions smoothly into a stable horizontal line at rest, completely devoid of spurious numerical oscillations. 

\begin{figure}[htbp]
    \centering
    \includegraphics[width=1.0\textwidth]{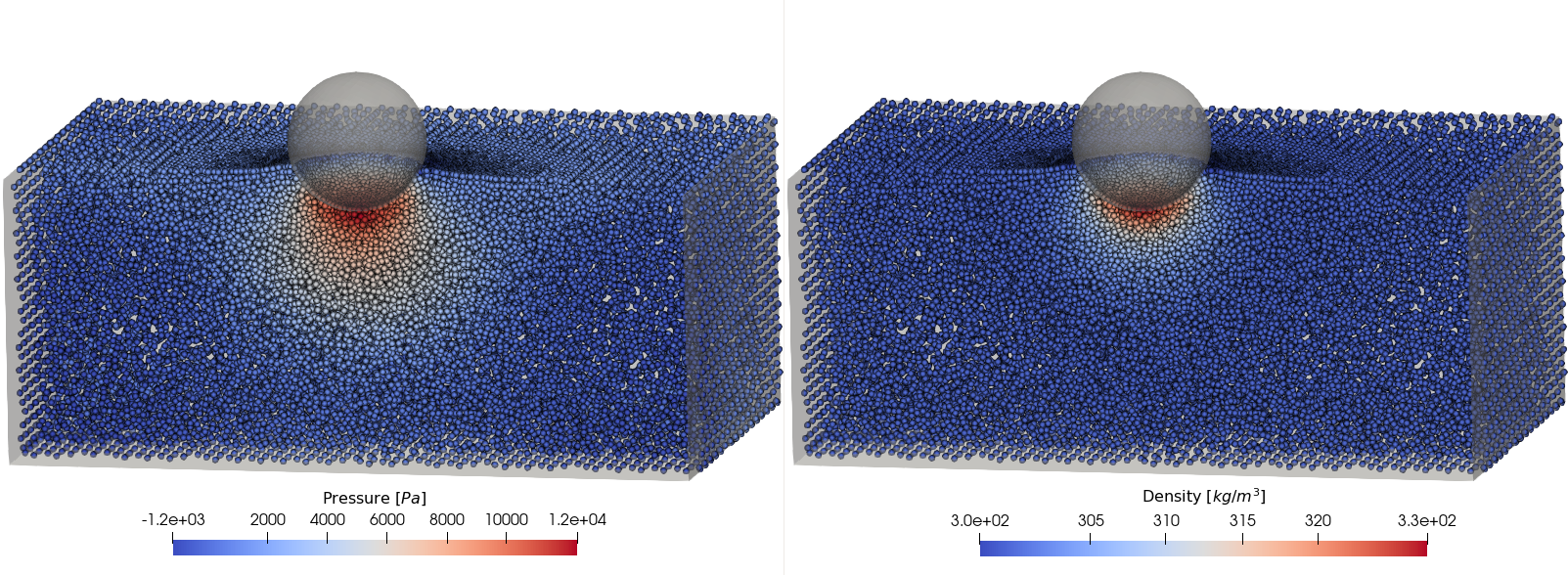} 
    \caption{Snow pressure and density field at time $t=0.105$ during the FSI simulation of a sphere falling into the snow.}
    \label{fig:spheredrop-density-pressure}
\end{figure}

Finally, Figure~\ref{fig:spheredrop-density-pressure} depicts the spatial distributions of the pressure and density fields at time $t=0.105$ during the snow compaction phase. Supplementary Video 3 shows the temporal evolution of these fields over the entire simulation duration, until the sphere comes to a halt. These contour plots clearly illustrate the generation of a globally smooth pressure field across the collocation points, validating the theoretical advantages of the proposed mixed formulation discussed in the previous section~\ref{sec:MCC_GFDM}.

\subsection{Snow Structure Interaction (Vehicle Application)}
\begin{figure}[htbp]
    \centering
    \includegraphics[width=0.9\textwidth]{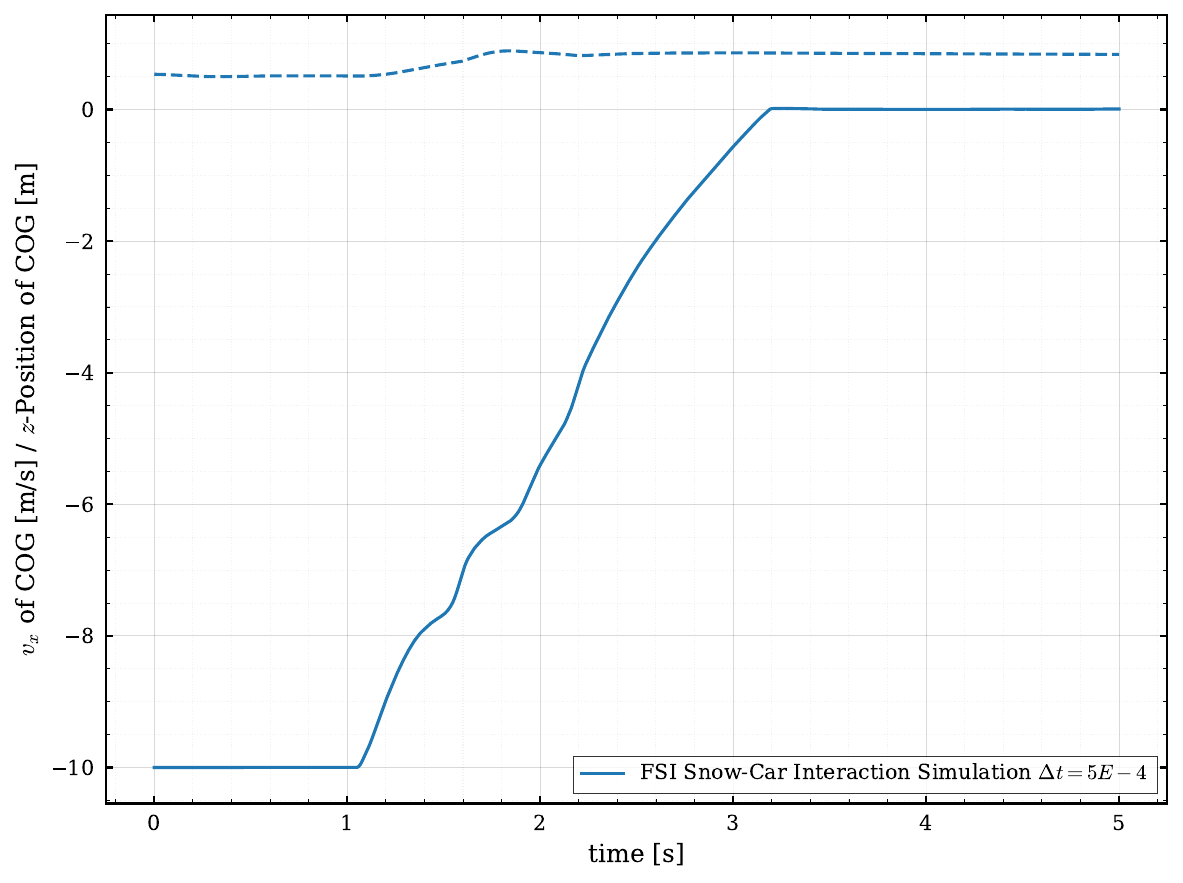} 
    \caption{The vertical displacement of the vehicle's center of gravity and its velocity in the direction of travel.}
    \label{fig:COGcar}
\end{figure}

\begin{figure}[htbp]
    \centering
    \includegraphics[width=0.8\textwidth]{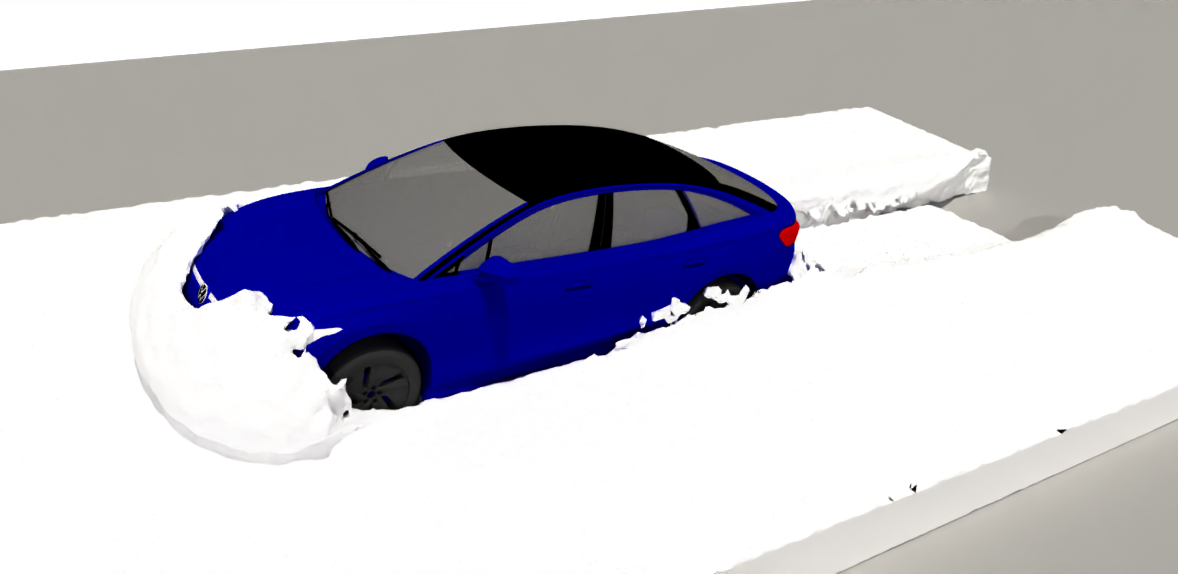} 
    \caption{FSI simulation of a car driving into a snowpack ($t=2$).}
    \label{fig:car}
\end{figure}

\begin{figure}[htbp]
    \centering
    \includegraphics[width=0.9\textwidth]{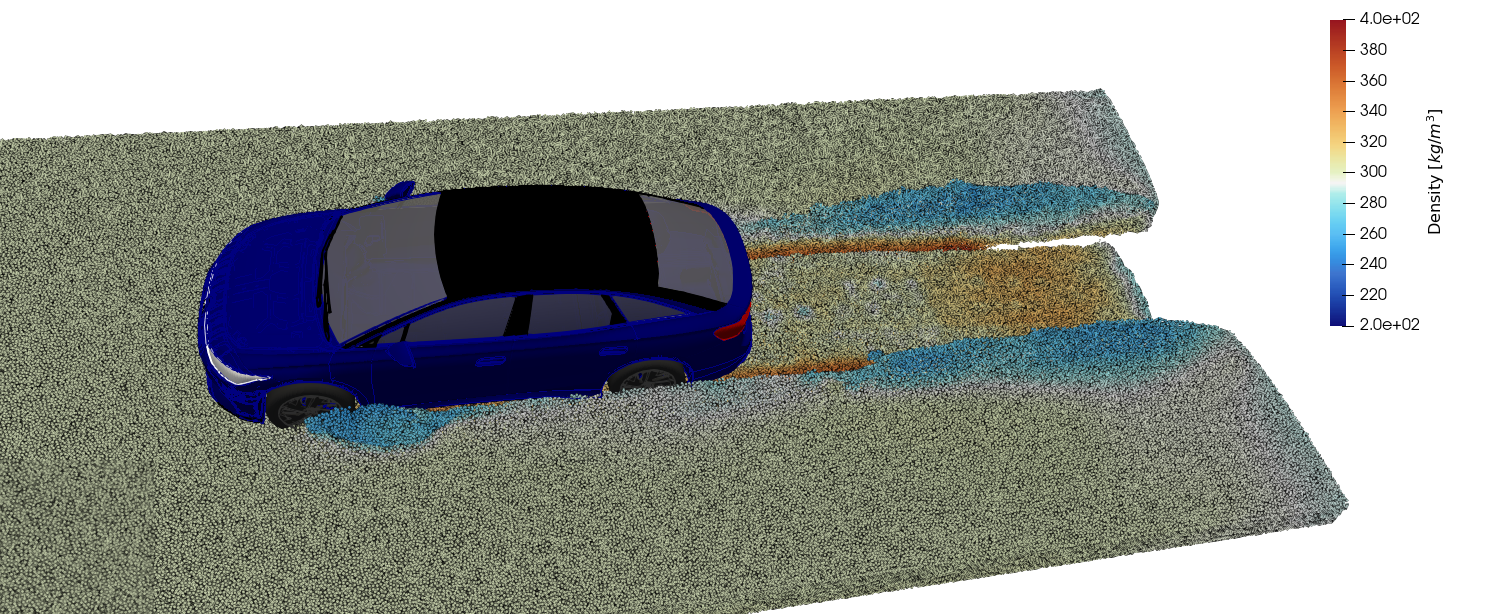} 
    \caption{Snow density field during the FSI simulation of a car driving into a snowpack ($t=4$). The field exhibits localized compaction within the tire tracks, leading to an increase in density, whereas shear-induced displacement by the vehicle body causes a density reduction adjacent to the tracks.}
    \label{fig:car-density}
\end{figure}

As a final demonstrative example, we present a purely qualitative simulation to highlight the potential of the coupled framework for full-scale automotive applications, which serves as the overarching motivation of this work. To this end, a complete vehicle, modeled as a rigid body, is propelled into a snowbank. The simulation captures the deceleration of the vehicle's initial velocity due to the snow-structure interaction forces. Furthermore, it demonstrates that, owing to the strain hardening and compaction intrinsic to the snow model, the vehicle is ultimately capable of riding up onto the consolidated snowpack. The trajectory of the vehicle's center of gravity and its velocity are depicted in Figure~\ref{fig:COGcar}. Supplementary Video 4 illustrates the vehicle driving into the snowpack, while Supplementary Video 5 displays the simultaneous density evolution within the snow. Corresponding snapshots at time $t=2$ and $t=4$ are shown in Figure~\ref{fig:car} respectively Figure~\ref{fig:car-density}. The stable execution of such a highly nonlinear scenario underscores the overall robustness and capability of the present coupled framework. Ultimately, this showcases the method's capacity for seamlessly handling complex contact cases and the intricate, large deformation physical interactions inherent to real world automotive environments.

    \section{Conclusions}\label{sec:conclusions}
This work presents the first successful integration of a complex, thermodynamically consistent constitutive model for snow into a purely meshless collocation framework. Diverging from conventional solvers that rely on the explicit integration of the full stress tensor, the proposed methodology employs a globally implicit calculation of the pressure field, which is then constitutively corrected using the return mapping scheme from \cite{Gaume2018} respectively \cite{Blatny2024}. Furthermore, the deviatoric stress is treated semi-implicitly by leveraging a numerical viscosity formulation. This decoupled numerical strategy is highly advantageous: It ensures the stability of significantly larger time steps and yields exceptionally smooth spatial pressure fields. Consequently, the accurate evaluation of boundary forces on complex geometries is greatly enhanced, thereby guaranteeing an optimal representation of the ensuing rigid body kinematics.

Moving forward, this framework provides a robust foundation for simulations of vehicle-snow interactions and enables the precise mapping of dynamic loads onto structural components. However, to fully utilize dynamic capabilities like adaptive spatial refinement, the implemented softening behavior necessitates proper regularization to circumvent localization artifacts. The macroscopic fracture representation can be further enriched by coupling the current formulation with a dedicated damage model (compare with \cite{Moeineddin2024}). Crucially, our decoupled viscosity formulation offers a flexible architecture. It enables the implicit integration of natural, physical viscosity models as an elegant alternative to enforcing an overstress formulation within the Cam-Clay return mapping procedure. Accordingly, a promising avenue for future research is the implementation of an inertia-dependent viscosity to rigorously resolve the inherent rate dependency of snow. By incorporating this physical formulation, the framework will guarantee a consistent macroscopic sinkage in section~\ref{sec:spheredrop}, effectively eliminating unphysical artificial dependencies caused by timestep refinement or increased penetration velocities.

\newpage

    \clearpage
    \section*{Declaration of Generative AI and AI-assisted technologies in the writing process}
    During the preparation of this work the authors used ChatGPT and Gemini in order to assist the writing. After using this tool, the authors reviewed and edited the content as needed and take full responsibility for the content of the published article.
    \section*{Declaration of Competing Interest}
    The authors declare that they have no known competing financial interests or personal relationships that could have appeared to influence the work reported in this paper.
    \section*{Acknowledgements}
    The authors would like to acknowledge Prof. Johan Gaume (ETH Zurich/ SLF Davos) for insightful comments and constructive discussions during the early stages of this work.

    \clearpage

    \bibliography{mybibfile.bib}

    \clearpage

    \appendix
\section{Derive Conservation Equations for Primitive Variables}\label{appendix:conservation-eq}
The conserved quantities are mass (represented by the density $\rho$), momentum ($\rho \mathbf{v}$), and total energy ($\rho E$). The governing equations formulated in Lagrangian framework are
\begin{align}
    \frac{D\rho}{Dt}&= - \rho \left(\nabla \cdot \mathbf{v}\right) , \nonumber\\
    \frac{D(\rho\mathbf{v})}{Dt} + \rho \mathbf{v} \left(\nabla\cdot\mathbf{v}\right) &= \nabla \cdot \boldsymbol{\sigma} + \rho\mathbf{g},\label{eq:continuum-energy} \\
     \frac{D(\rho E)}{Dt} + \rho E \left(\nabla\cdot\mathbf{v}\right) &= \nabla \cdot (\boldsymbol{\sigma}\mathbf{v}) + \rho\mathbf{g}\cdot\mathbf{v}, \nonumber
\end{align}
We rewrite this set of equations in primitive variables, i.e. density, velocity and pressure, together with the deviatoric stress tensor. To achieve this, we replace the equation for the total energy $\rho E$ by an equation for the Cauchy stress. We note that the total energy splits into potential and kinetic energy  
\begin{align}
    \rho E =  \rho e + \frac{1}{2}\rho \mathbf{v}^2. \label{eq:total-energy}
\end{align}
where $e$ is the specific energy (or the rate of mechanical work per unit mass). We plug the rate form of equation~\eqref{eq:total-energy},
\begin{align}
\frac{D(\rho E)}{Dt} = \frac{(D\rho e)}{Dt} + \mathbf{v}\frac{(D\rho\mathbf{v})}{Dt}    
\end{align}
as well as the term
\begin{align}
\frac{D(\rho e)}{Dt} = \rho\frac{De}{Dt} + e\frac{D\rho}{Dt}.
\end{align}
into the last equation of system~\eqref{eq:continuum-energy}. Further sorting of terms, using mass and momentum conservation yields the energy balance equation as the rate of specific energy $e$    
\begin{align}
   \frac{De}{Dt}= \frac{1}{\rho}\left(\boldsymbol{\sigma} : \nabla\mathbf{v}\right). \label{eq:inner-energy}
\end{align}
In general there could be a dissipation term $\mathcal{D}$ on the right hand side of equation~\eqref{eq:inner-energy} subjected to the Clausius-Duhem inequality, which represents the continuum mechanics formulation of the second law of thermodynamics. This inequality strictly requires that the dissipated energy remains non-negative ($\mathcal{D} \geq 0$). In the numerical solution of finite strain plasticity, adherence to this thermodynamic restriction is mathematically governed by the Kuhn-Tucker loading conditions, which will be defined by the yield function in section~\ref{subsec:cam-clay}. These conditions explicitly control the transition between pure elasticity ($\mathcal{D} = 0$) and active plastic yielding, which drives the irreversible dissipation. Consequently, the Kuhn-Tucker conditions ensure that the numerical integration of the plastic flow remains thermodynamically consistent at all times (\cite{Simo1998}, \cite{Holzapfel2000}).

\section{Thermodynamically Consistent Derivation of the Constitutive Model}\label{appendix:thermodynamic-consistency}

In this appendix, we formally demonstrate that the chosen strain energy density function leads to a thermodynamically consistent hyperelastic constitutive relation. This also induces the objective rate formulations used in section~\ref{subsec:hypoelastic} and density-dependent global moduli for the Cauchy stress used in section~\ref{subsubsec:elastic-trial}.

We begin by recalling the specific strain energy per unit mass $e$, formulated in terms of the logarithmic Hencky strain $\varepsilon^E$:
\begin{equation}
     e(\eps^E) := \frac{1}{\rho_0}\left(\frac{1}{2}\lambda_{\tau} \tr(\eps^E)^2 + \mu_{\tau} \tr\left((\eps^E)^2\right)\right).\label{appendix-eq:e(eps)}
\end{equation}
Assuming pure elastic forces, the local energy balance equation in terms of the Cauchy stress $\bsigma$ and the rate of deformation tensor $\mathbf{D}^E$ is given by:
\begin{equation}
    \frac{De}{Dt} = \frac{1}{\rho}\left(\boldsymbol{\sigma} : \mathbf{D}^E\right). 
\end{equation}
Multiplying both sides by the reference density $\rho_0$ allows us to rewrite the right-hand side in terms of the Kirchhoff stress $\boldsymbol{\tau} := \frac{\rho_0}{\rho}\bsigma = J\bsigma$
\begin{equation}
    \rho_0 \frac{De}{Dt} = \boldsymbol{\tau} : \mathbf{D}^E.
\end{equation}
Utilizing the objective logarithmic rate $ \overset{\circ}{{\eps}}^{E,\log} =\mathbf{D}^E$, we obtain
\begin{equation}
    \rho_0 \frac{De}{Dt} = \frac{\partial (\rho_0 e)}{\partial \eps^E} : \overset{\circ}{{\eps}}^{E,\log} = \frac{\partial (\rho_0 e)}{\partial \eps^E} : \mathbf{D}^E.
\end{equation}
Equating these two expressions for the specific energy rate yields the general hyperelastic constitutive relation
\begin{equation}
    \boldsymbol{\tau} = \frac{\partial (\rho_0 e)}{\partial \eps^E}.
\end{equation}
To find the explicit stress response, we perform the tensor derivative of the strain energy potential $\Psi = \rho_0 e$, which leads to
\begin{equation}
    \frac{\partial (\rho_0 e)}{\partial \eps^E} = \frac{1}{2} \lambda_\tau \cdot 2 \tr(\eps^E) \frac{\partial \tr(\eps^E)}{\partial \eps^E} + \mu_\tau \frac{\partial \tr((\eps^E)^2)}{\partial \eps^E}.
\end{equation}
Using the standard tensor derivative identities $\frac{\partial \tr(\eps^E)}{\partial \eps^E} = \mathbf{I}$ and $\frac{\partial \tr((\eps^E)^2)}{\partial \eps^E} = 2\eps^E$, this evaluates to the hyperelastic law for the Kirchhoff stress
\begin{equation}
    \boldsymbol{\tau} = \lambda_\tau \tr(\eps^E)\mathbf{I} + 2\mu_\tau \eps^E.
\end{equation}
This derivation confirms that the formulation is consistent with the second law of thermodynamics.

In order to obtain the formulation for the actual physical Cauchy stress $\bsigma$, we pull back the Kirchhoff stress using the volume ratio $J = \frac{\rho_0}{\rho}$:
\begin{equation}
    \bsigma = J^{-1}\boldsymbol{\tau} = \lambda \tr(\eps^E)I + 2\mu \eps^E
\end{equation}
This explicitly defines the density-dependent Lamé parameters $\lambda = J^{-1}\lambda_\tau$ and the shear modulus $\mu = J^{-1}\mu_\tau$. Based on these parameters, we can derive the density dependent bulk modulus $K = \lambda + \frac{2}{3}\mu = J^{-1}K_\tau$.

To relate this hyperelastic model to a hypoelastic formulation for numerical implementation, we take the objective logarithmic rate of the Cauchy stress. Accounting for the volume change rate $\dot{J} = J \tr(\mathbf{D}^E)$, we obtain
\begin{equation}
    \overset{\circ}{\boldsymbol{\sigma}}^{\log} =  \lambda\tr(\mathbf{D}^E)\mathbf{I} + 2\mu\mathbf{D}^E - \tr(\mathbf{D}^E)\boldsymbol{\sigma} .\label{appendix-eq:hypoelasticlaw-cauchy}
\end{equation}
Converting this to the material time derivative of the Cauchy stress yields
\begin{equation}
    \dot{\boldsymbol{\sigma}} =  \lambda\tr(\mathbf{D}^E)\mathbf{I} + 2\mu\mathbf{D}^E -\bsigma\boldsymbol{\Omega}^{\log}+\boldsymbol{\Omega}^{\log}\bsigma- \tr(\mathbf{D}^E)\boldsymbol{\sigma}.
\end{equation}
To isolate the volumetric response, we take $-1/3$ times the trace on both sides, providing the rate of the hydrostatic pressure $p$
\begin{equation}
    \dot{p} = -(K + p) \tr(\mathbf{D}^E).
\end{equation}
Finally, utilizing the mass conservation equation $\dot{\rho} + \rho \tr(\mathbf{D}^E) = 0$, we can substitute the trace of the rate of deformation tensor to rearrange the pressure rate into an elastic equation of state. This establishes the analytic elastic compressibility modulus
\begin{equation}
    \kappa_E = \frac{1}{K + p} = \frac{1}{J^{-1}K_\tau + p}.
\end{equation}

\end{document}